\documentclass[a4paper,12pt,reqno]{amsart}
\usepackage{amssymb}
\usepackage{delarray}

\newtheorem{thm}{Theorem}[section]
\newtheorem{cor}[thm]{Corollary}
\newtheorem{lem}[thm]{Lemma}
\newtheorem{prop}[thm]{Proposition}
\theoremstyle{definition}

\theoremstyle{remark}
\newtheorem{rem}[thm]{Remark}
\numberwithin{equation}{section}
\newtheorem{example}{Example}[section]
\newcommand{\Real}{\mathbb R}

\newcommand{\eqnref}[1]{(\ref{#1})}

\newcommand{\V}{{\bf v}}

\newcommand{\ms}{{\mathfrak{sl}_2}}

\def\var{\varphi}

\def\l{\lambda}

\def\ep{\epsilon}

\newcommand{\px}{\partial_x}

\newcommand{\pt}{\partial_t}
\newcommand{\pu}{\partial_u}

\newcommand{\tx}{\tilde{x}}
\newcommand{\tilt}{\tilde{t}}

\begin{document}

\title[New Integral Equations]{New Classes of Non-Convolution Integral Equations Arising From Lie Symmetry Analysis of Hyperbolic PDEs}%
\author{Mark Craddock}
\address{School of Mathematical Sciences\\ University of Technology Sydney\\ PO Box 123, Broadway\\ New South
Wales 2007\\ Australia}
\email{Mark.Craddock@uts.edu.au}
\author{Semyon Yakubovich}
\address{Department of Pure Mathematics, Faculty of Sciences, University of Porto, Campo Alegre
St., 687, 4169-00,   Porto, Portugal.}

\email{syakubov@fc.up.pt}
\date{September 2015}%


\begin{abstract}
In this paper we consider some new classes of integral equations
that arise from Lie symmetry analysis. Specifically, we consider the
task of obtaining solutions of a Cauchy problem for some classes of
second order hyperbolic partial differential equations. Our analysis
leads to new integral equations of non-convolution type, which can
be solved by classical methods. We derive solutions of these
integral equations, which in turn lead to solutions of the
associated Cauchy problems.
\end{abstract}
\maketitle

\section{Introduction}

The theory of continuous symmetry groups of systems of partial
differential equations was developed by Lie in the last decades of
the nineteenth century. A symmetry is a transformation which maps
solutions to other solutions. Lie developed a method for
systematically computing all continuous symmetries of a given system
of differential equations. Excellent modern accounts may be found in
the books by Olver \cite{Olv93} and those of Bluman and his
coauthors, such as the volume \cite{BK89}.

Symmetries are powerful tools and they allow us to solve a wide
variety of problems. For example, we may compute fundamental
solutions of parabolic equations from a knowledge of their
symmetries. See the papers \cite{CL2006},\cite{CL2007},
\cite{Cra2008}, \cite{CL2012} for details, together with many
examples and applications. Boundary value problems may also be
solved by the method of group invariant solutions. The
aforementioned book \cite{BK89} contains an extensive discussion of
this topic. Symmetries are also essential to the study of
conservation laws. A chapter of \cite{Olv93} is devoted to this.

A new way of applying Lie symmetries was  introduced in
\cite{Cra2008} and developed extensively in \cite{Cra2014}. The idea
is to construct an integral operator that maps test functions to
solutions of the PDE by integration against a symmetry solution. The
purpose of this paper is to show how the method can lead to new
classes of non-convolution integral equations, which may be solved
by a combination of classical techniques.

We first provide an example to illustrate how the method works. We
then turn to a class of second order linear hyperbolic equations and
derive some new integral equations which arise in the solution of
Cauchy problems associated to these equations. We then solve these
equations using a novel combination of integral transform methods.

The outline of the paper is as follows. In Section 2 we introduce
the method that forms the basis of our analysis. In Section 3 we
determine the symmetries of a class of hyperbolic PDEs of the form
$u_{tt}=u_{xx}+f(x)u$. We show that there are nontrivial symmetries
when $f$ satisfies one of three families of Riccati equations. In
Section 4 we study the equation $u_{tt}=y^2u_{yy}+yu_y+y^2u$ ( which
is equivalent to $u_{tt}=u_{xx}+e^{2x}u$, with $y=e^x$) and show how
the symmetries lead to a solution of the Cauchy problem for this
equation via new non-convolution integral equations, which we are
able so solve. See Theorems 4.2 and 4.3. In Section 5 we set up the
integrals for the Cauchy problem for the more general equation
$u_{tt}=u_{xx}+(Ae^x+Be^{-x})^{-2}u$. In Section 6 we solve the
integral equation in the case $A=B=1/2$ and obtain a solution of the
problem $u_{tt}=u_{xx}+(\mathrm{sech}^2x)u, x>0$, subject to
$u(x,0)=f(x),
 u_t(x,0)=0.$ (See Theorem 6.6). A brief conclusion follows. We believe that the
methods of this paper may prove to be of use in the study of various
types of linear PDEs.

\section{Introduction to the Method}
We illustrate the approach by solving the problem
\begin{align*}
u_{tt}= u_{xx}-\frac{A}{x^2}u,\ u(x,0)=f(x), u_t(x,0)=g(x),\ x> 0,
t\geq0.
\end{align*}
This problem was solved using different symmetries in
\cite{Cra2014}. Here we use a simpler approach. We suppose that
$f,g\in \mathcal{S}(\Real_+)$ the Schwartz space, consisting of
smooth functions whose derivatives (including the function itself)
decay at infinity faster that any power, i.e. Schwartz functions are
rapidly decreasing.  It is known that the Fourier Transform and
Mellin transforms are automorphisms of the Schwartz space, see
\cite{SW71}. This space is a topological vector space of functions
$\varphi$ such that $ \varphi \in C^\infty (\mathbb{R}_+)$ and
$x^\alpha \varphi^{(\beta)} (x) \to 0, \ x \to \infty,\ \alpha,
\beta \in \mathbb{N}_0$.

It is elementary that the scaling transformation $u(x,t)\to u(\l
x,\l t)$ is a symmetry. The idea is to introduce a new solution by
setting
\begin{align}
U(x,t)=\int_0^\infty \varphi(\l)u(\l x,\l t)d\l.
\end{align}

With separation constant equal to one half, separation of variables
leads to the solution $u_1(x,t)=\sqrt{x}I_\nu(x)e^{t},$ where
$\nu=\frac{1}{2} \sqrt{4 A+1}$ and $I_\nu(x)$ denotes the usual
modified Bessel function, (see 9.6.10  of \cite{AS72}). We denote
the Bessel function of the first kind by $J_\nu(z)$, (9.1.10 of
\cite{AS72}).  Applying the symmetry and replacing $\l$ with $i\l$,
we obtain a solution
\begin{align*}
u(x,t)&=\int_0^\infty\varphi(\l)\sqrt{\l x}J_\nu(\l x)\cos(\l
t)d\l\\& \quad +\int_0^\infty\psi(\l)\sqrt{\l x}J_\nu(\l x)\sin(\l
t)d\l.
\end{align*}
We absorbed the terms in $i$ into $\varphi$ and $\psi$. The
requirement that $u(x,0)=f(x)$ leads to the equation
\begin{align*}
f(x)=\int_0^\infty\varphi(\l)\sqrt{\l x}J_\nu(\l x)d\l.
\end{align*}
We immediately recognise that $f$ must be the Hankel transform of
$\varphi$. The inversion theorem for the Hankel transform, (see
\cite{tit}) yields

\begin{align*}
\varphi(\l)=\int_0^\infty f(y)\sqrt{\l y}J_\nu(\l y)dy.
\end{align*}
For the condition $u_t(x,0)=g(x)$ we find
\begin{align}
\int_0^\infty \l\psi(\l)\sqrt{\l x}J_\nu(\l x)d\l=g(x),
\end{align}
which gives
\begin{align}
\psi(\l)=\frac1\l\int_0^\infty g(y)\sqrt{\l y}J_\nu(\l y)dy.
\end{align}

Hence

\begin{align*}
u(x,t)&= \int_0^\infty \int_0^\infty f(y)\sqrt{xy}\l J_\nu(\l
x)J_\nu(\l y)\cos(\l t)dy d\l\\&+\int_0^\infty \int_0^\infty
g(y)\sqrt{xy}  J_\nu(\l x)J_\nu(\l y)\sin(\l t)dy d\l,
\end{align*}
solves our initial value problem.  Thus we have obtained a solution
of the initial value problem from a separable solution and a scaling
symmetry. This solution can be simplified by making use of the
result
\begin{align*}
\int_0^\infty &\sin(cx) J_{\mu}(ax)J_\mu(bx)dx\\   & = 0, \hskip5.8cm\ 0<c<b-a, 0<a<b\\
&
=\frac{1}{2\sqrt{ab}}P_{\mu-1/2}\left(\frac{b^2+a^2-c^2}{2ab}\right),\
\ \
b-a<c<b+a, 0<a<b\\
&=\frac{\cos(\mu\pi)}{\pi\sqrt{ab}}Q_{\mu-1/2}\left(-\frac{b^2+a^2-c^2}{2ab}\right),\
b+a<c,0<a<b,
\end{align*}
which is valid for $\mu>-1$. (See formula 6.672.1 of \cite{GR2000}).
Here $P_{\nu}$ and $Q_\nu$ are Legendre functions, \cite{OLBC2010}.
We define
\begin{align}
p(x,y,t)=\int_0^\infty J_\nu(\l x)J_\nu(\l y)\sin(\l t) d\l
\end{align}
via the preceding integral and the solution can be written
\begin{align}\label{sao}
u(x,t)&= \int_0^\infty \int_0^\infty f(y)\sqrt{xy}\l J_\nu(\l
x)J_\nu(\l y)\cos(\l t)dy d\l\nonumber\\&+ \int_0^\infty
\sqrt{xy}g(y)p(x,y,t)dy.
\end{align}
The absolute convergence of the second integral justifies our use of
Fubini's Theorem.

The integral $\int_0^\infty  \l J_\nu(\l x)J_\nu(\l y)\cos(\l t)
d\l$ is divergent, so we cannot reverse the order of integration in
the first integral in \eqnref{sao}. However, we can define a
solution in the distributional sense by  noting that for each $n$
$$\int_0^n \l J_\nu(\l x)J_\nu(\l y)\cos(\l t)
d\l= \frac{\partial}{\partial t}\int_0^n   J_\nu(\l x)J_\nu(\l
y)\sin(\l t)d\l.$$ Taking the limit as $n\to\infty$ we define the
distribution $\Lambda$ on $\mathcal{S}(\Real_+)$ by setting
\begin{align}
\Lambda(f)=\int_0^\infty f(y)\frac{\partial}{\partial
t}\int_0^\infty J_\nu(\l x)J_\nu(\l y)\sin(\l t)d\l dy.
\end{align}
This leads to a solution in the distibutional sense given by
\begin{align}
u(x,t)&=  \frac{\partial}{\partial t} \int_0^\infty
f(y)\sqrt{xy}p(x,y,t)dy  + \int_0^\infty \sqrt{xy}g(y)p(x,y,t)dy.
\end{align}

In general this procedure requires us to solve one or more integral
equations if we are to satisfy the initial data. Often this equation
turns out to be a familiar integral transform with a known inversion
formula, as happened with our example. In this paper we will show
how it can lead to   entirely new solvable integral equations.

This method hints at a deep connection between Lie symmetries and
classical harmonic analysis and indeed this is the case. The
relationship between Lie symmetry analysis and harmonic analysis and
representation theory has been developed over the years, beginning
with \cite{Cra95}, \cite{Cra94} and \cite{Cra2000}. See also
\cite{CD2009}, \cite{CL2012}, \cite{SS2010}, \cite{SS2004} and
\cite{Fra2012}.

The basic idea of this work is that for linear PDEs, it is often
possible to realise the Lie symmetries as global representations of
the underlying Lie group in the following sense. One has a Lie group
$\mathcal{G}$, a representation $(\sigma,V)$ of $\mathcal{G}$ and a
mapping $A:V\to H$ where $H$ is a solution space of the equation. So
$A$ maps $f\in V$ to some solution of the PDE. Then if the action of
the symmetries on solutions $u=Af$, is denoted by $\rho(g)u$ for
$g\in \mathcal{G},$ the following relationship holds.

\begin{align}
(\rho(g)Af)(x)=(A\sigma(g)f)(x),\ g\in \mathcal{G},
\end{align}
where $f\in V$

In fact this relationship has been established for many classes of
PDEs and the previous references detail some of the work done in
this area.

A consequence of all this is that the operator
\begin{align*}
v_i(x)&=\int_\Omega\var(\ep)\rho(\exp\ep\V_i)u(x)d\ep,\nonumber\\
\end{align*}
is essentially a group theoretic Fourier transform, which inherits
properties of the representations equivalent to $\rho$. It is
therefore not surprising that this technique so often leads to the
problem of inverting a standard integral transform. In this paper we
consider some problems where the integral equations do not appear in
the existing literature.

\section{The Equation $u_{tt}=u_{xx}+f(x)u$}

We will determine the forms of the potential function $f$ which
allow non trivial symmetries. For an in depth discussion of the
computation of Lie symmetries, see the book by Olver \cite{Olv93}.
We look for a vector field
$\V=\xi(x,t)\px+\tau(x,t)\pt+\phi(x,t,u)\pu $ which generates the
Lie symmetries of the PDE. An application of Lie's Theorem (Theorem
2.31 p104, \cite{Olv93}) leads to the defining equations
\begin{align}\label{stex2}
\phi^{tt}=\phi^{xx}+f(x)\phi+f'(x)\xi u,
\end{align}
where
\begin{align*}
\phi^{xx}&=\phi_{xx}+(2\phi_{xu}-\xi_{xx})u_x-\tau_{xx}u_t+\phi_{uu}u_x^2+(\phi_u-2\xi_x)
u_{xx}-2\tau_xu_{xt}
\\
\phi^{tt}&=\phi_{tt}-\xi_{tt}u_x+(2\phi_{tu}-\tau_{tt})u_t+\phi_{uu}u_t^2-2\xi_tu_{xt}+(\phi_u-2\tau_t)u_{tt}.
\end{align*}

From the defining equations \eqnref{stex2}, we find the conditions
\begin{align*} \phi&=\alpha u+\beta,\ \xi_t=\tau_x,
-\xi_{tt}=2\phi_{xu}-\xi_{xx},\\
2\phi_{ut}&=\tau_{tt}-\tau_{xx}, \phi_{u}-2\tau_t=\phi_u-2\xi_x\\
\phi_{tt}&+(\phi_u-2\tau_t)f(x)u=\phi_{xx}+f(x)\phi+f'(x)u\xi.
\end{align*}
The function $\beta$ will be an arbitrary solution of the PDE.
Nothing further can be said about it. From the first and last of
these equations we have
\begin{align}\label{khawes}
\alpha_{tt}=\alpha_{xx}+f'(x)\xi+2\tau_tf(x).
\end{align}
It is also easy to see that $\xi_{xx}=\xi_{tt}$ and
$\tau_{tt}=\tau_{xx}.$ Consequently $2\phi_{ut}=0$, and this means
$\alpha$ is independent of $t.$ Similarly we see that $\alpha$ is
also independent of $x.$ It is thus a constant. From \eqnref{khawes}
we see that
\begin{align}
2\xi_x+\frac{f'}{f}\xi=0.
\end{align}
This means that $\xi(x,t)=e^{-\frac12\ln
f}c(t)=\frac{c(t)}{\sqrt{f}}$, where $c$ has yet to be fixed. With
$G(x)=\frac{c(t)}{\sqrt{f}}$, we now require
\begin{align}
c''(t)G(x)=G''(x)c(t).
\end{align}
This is only possible if $G''/G=\l$ a constant. We then have three
possibilities for $c:$
\begin{align}
c(t)=\begin{cases}&c_1t+c_2\\
&c_1\sin(\omega t)+c_2\cos(\omega t)\\
& c_1e^{\omega t}+c_2e^{-\omega t}.
\end{cases}
\end{align}
Clearly the three possible choices for $G$ are
\begin{align}
G(x)=\begin{cases}&Ax+B\\& Ae^{\omega x}+Be^{-\omega x}\\ &
A\cos\omega x+B\sin\omega x
\end{cases}
\end{align}
where $A,B$ are arbitrary constants. Correspondingly, the possible
choices for $f$ are
\begin{align}
f(x)=\begin{cases}& \frac{1}{(Ax+B)^2}\\&\frac{1}{(Ae^{\omega
x}+Be^{-\omega x})^2}\\ &\frac{1}{(A\cos \omega x+B\sin \omega
x)^2}.
\end{cases}
\end{align}

We can now write down a basis for the Lie algebra of symmetries for
each choice of $f$.  In each case, the Lie algebra is isomorphic to
$\ms\oplus\mathfrak{r}$, where $\mathfrak{r}$ is the Lie algebra
generating the additive Lie group $\Real.$

We are interested in the second case. The first case is obtainable
from our illustrative example by a linear change of variables. If
$$f(x)=\frac{1}{(Ae^{\omega x}+Be^{-\omega x})^2},$$ a basis for the
Lie symmetry algebra is
\begin{align*}
\V_1&=e^{\omega t}(A e^{\omega x}+Be^{-\omega x})\px+e^{\omega
t}(Ae^{\omega x}-B e^{-\omega x})\pt,\\ \V_2&=e^{-\omega t}(A
e^{\omega x}+Be^{-\omega x})\px-e^{-\omega t}(Ae^{\omega x}-B
e^{-\omega x})\pt,\\ \V_3&=\pt, \V_4=u\pu,\V_{\beta}=\beta\pu.
\end{align*}

The vector field $\V_3$ clearly induces the symmetry $u(x,t)\to
u(x,t+\ep).$ Let us obtain the symmetry arising from the first
vector field. We have to solve
\begin{align}
\frac{d\tx}{d\ep}&=e^{\omega\tilt}(Ae^{\omega\tx}+Be^{-\omega\tx}),\ \tx(0)=x,\\
\frac{d\tilt}{d\ep}&=e^{\omega\tilt}(Ae^{\omega\tx}-Be^{-\omega\tx}),\
\tilt(0)=t.
\end{align}
Adding these equations gives $$
\frac{d(\tx+\tilt)}{d\ep}=2Ae^{\omega(\tx+\tilt)}. $$

Now setting $Z=\tx+\tilt$, we solve this first order DE and get
\begin{align}\label{expsymz}
e^{\omega(\tx+\tilt)}=\frac{e^{\omega(x+t)}}{1-2A\ep\omega
e^{\omega(x+t)}}.
\end{align}

Next we solve
\begin{align}
\frac{d(\tx-\tilt)}{d\ep}=2B e^{-\omega(\tx-\tilt)}.
\end{align}
With $z=\tx-\tilt$, $z(0)=x-t$ we obtain the expression
\begin{align}
e^{\omega(\tx-\tilt)}=e^{\omega(x-t)}+2B\ep\omega.
\end{align}
Taking logs we have a pair of simultaneous equations for $\tx$ and
$\tilt.$ This gives us

\begin{align}
e^{\omega\tx}=\left(\frac{e^{2\omega x}+2B\ep\omega
e^{\omega(x+t)}}{1-2A\ep\omega e^{\omega(x+t)}}\right)^{\frac12}
\end{align}
and
\begin{align*}
e^{\omega\tilt}&=\frac{e^{\omega(x+t)}}{\sqrt{(e^{2\omega
x}+2B\ep\omega e^{\omega(x+t)})(1-2A\ep\omega e^{\omega(x+t)}})}.
\end{align*}

From these we obtain expressions for $\tilde{x}$ and $\tilde{t}.$
From this we conclude that if $u$ is a solution of
\begin{align}
u_{tt}=u_{xx}+\frac{1}{(Ae^{\omega x}+Be^{-\omega x})^2}u,
\end{align}
then so is
\begin{align}\label{symmetry}
\tilde{u}(x,t;\ep)=&u\left( \frac12\ln\left(\frac{e^{2\omega
x}+2B\ep\omega e^{\omega(x+t)}}{1-2A\ep\omega
e^{\omega(x+t)}}\right), \right.\nonumber\\&\left.
\ln\left(\frac{e^{\omega(x+t)}}{\sqrt{(e^{2\omega x}+2B\ep\omega
e^{\omega(x+t)})(1-2A\ep\omega e^{\omega(x+t)}})}\right)\right).
\end{align}

The vector field $\V_2$ can be exponentiated by similar means. This
is left to the interested reader.
\section{The Case of one zero constant}

To begin our analysis we will first consider the special case when
$A=0,B=1$. The case $B$ arbitrary can be obtained from this by an
obvious change of variables. This leads to the equation
\begin{align}
u_{tt}=u_{xx}+e^{2x}u,\ \ x\in\Real.
\end{align}
For convenience we will make the change of variables $y=e^{x}.$ This
leads to the equation
\begin{align}\label{PDE1}
u_{tt}=y^2u_{yy}+yu_y+y^2u,\ \ y>0.
\end{align}
We will solve this subject to the initial conditions $u(y,0)=f(x),$
and $ u_t(y,0)=g(y)$ where $f,g$ are suitable functions. It is
straightforward to see that a stationary solution is
$u_0(y)=J_0(y)$. Application of the symmetry \eqnref{symmetry} under
the change of variables shows that
\begin{align}
\tilde{u}_1(y,t;\l)=J_0(\sqrt{y^2+2\l ye^t}),
\end{align}
is also a solution. Since $t\to -t$ is a symmetry, it is clear that
\begin{align}
\tilde{u}_2(y,t;\l)=J_0(\sqrt{y^2+2\l ye^{-t}})
\end{align}
is again a solution. We now consider
\begin{align}
u(y,t)=\frac12\int_0^\infty \varphi(\l)\left[J_0(\sqrt{y^2+2\l
ye^t})+J_0(\sqrt{y^2+2\l ye^{-t}})\right] d\l,
\end{align}
where $\varphi$ has suitable decay to guarantee convergence of the
integral. This is easily seen to satisfy the PDE \eqnref{PDE1} and
the conditions
\begin{align}
u(y,0)=\int_0^\infty \varphi(\l)J_0(\sqrt{y^2+2\l y})d\l
\end{align}
and $u_t(y,0)=0.$ Our first task is then to solve the integral
equation
\begin{align}\label{eqn3one}
\int_0^\infty \varphi(\l)J_0(\sqrt{y^2+2\l y})d\l=f(y).
\end{align}

We proceed by taking the Laplace transform in $y$ of both sides. We
suppose that $\varphi\in L_1(\Real_+)$. Via the known inequality for
the Bessel function $ \sqrt x |J_\nu(x) | < C,\ x>0$, where $C$ is
an absolutely positive constant, we have from \eqnref{eqn3one} the
estimate
\begin{align}\label{ineq}
|f(y)| \le C \int_0^\infty \frac{|\varphi(\lambda)| }{ (y^2+
2\lambda y)^{1/4} }d\lambda \le {C\over \sqrt y}
\|\varphi\|_{L_1(\mathbb{R}_+)}.
\end{align}
Thus
\begin{align*}
\left|\int_0^\infty f(y)e^{-sy}dy\right|\leq
C\|\varphi\|_{L_1(\Real_+)}\int_0^\infty
\frac{e^{-sy}}{\sqrt{y}}dy=\frac{\sqrt{\pi}C\|\varphi\|_{L_1(\Real_+)}}{\sqrt{s}}.
\end{align*}

Now 10.2.6, p59 of \cite{RK66} gives
\begin{align}
\int_0^\infty J_0(\sqrt{y^2+2\l
y})e^{-sy}dy=\frac{e^{-\l(\sqrt{s^2+1}-s)}}{\sqrt{s^2+1}}.
\end{align}
If $F(s) =\int_0^\infty f(y)e^{-sy}dy$ then an application of the
Laplace transform and Fubini's Theorem, (justified by
\eqnref{ineq}), gives
\begin{align}
\int_0^\infty
\varphi(\l)\frac{e^{-\l(\sqrt{s^2+1}-s)}}{\sqrt{s^2+1}}d\l =F(s) .
\end{align}
Letting $\Phi$ denote the Laplace transform of $\varphi$, we have
\begin{align}
\frac{\Phi(\sqrt{s^2+1}-s) }{\sqrt{s^2+1}}=F(s) .
\end{align}
The substitution $s=\sinh k$ reduces this to
\begin{align}
\Phi(e^{-k})=\cosh k F(\sinh k).
\end{align}
Putting $k=-\ln\xi$ gives
\begin{align}
\Phi(\xi)=\frac12\left(\xi+\frac1\xi\right)F\left(\frac12\left(\frac1\xi-
\xi\right)\right).
\end{align}
Thus
\begin{align}
\varphi(\l)=\mathcal{L}_\xi^{-1}\left[\frac12\left(\xi+\frac1\xi\right)F\left(\frac12\left(\frac1\xi-
\xi\right)\right)\right].
\end{align}

We have thus established the following.
\begin{prop}
The equation $u_{tt}=y^2u_{yy}+yu_y+y^2u$ with $u(y,0)=f(y)$ and
$u_t(y,0)=0$ has a solution
\begin{align}
u(y,t)=\frac12\int_0^\infty&\mathcal{L}_\xi^{-1}\left[\frac12\left(\xi+\frac1\xi\right)F\left(\frac12\left(\frac1\xi-
\xi\right)\right)\right](\l)\nonumber \\&\times
\left[J_0(\sqrt{y^2+2\l ye^t})+J_0(\sqrt{y^2+2\l ye^{-t}})\right]
d\l,
\end{align}
provided the integral converges. Here $F$ is the Laplace transform
of $f.$
\end{prop}

Our next task is to determine an expression for the inverse Laplace
transform. An elementary approach is to assume that $f$ is analytic.
That is
\begin{align}
f(y)=\sum_{n=0}^\infty \frac{1}{n!}f^{(n)}(0)y^n.
\end{align}
In this case
\begin{align}
F(s) =\sum_{n=0}^\infty  f^{(n)}(0)\frac{1}{s^{n+1}},
\end{align}
so that
\begin{align}
\frac12\left(\xi+\frac1\xi\right)F\left(\frac12\left(\frac1\xi-
\xi\right)\right)&=\sum_{n=0}^\infty f^{(n)}(0)
\frac{(\xi^2+1)(2\xi)^n}{(1-\xi^2)^{n+1}}.
\end{align}
Thus
\begin{align}
\varphi(\l)=\sum_{n=0}^\infty f^{(n)}(0)g_n(\l),
\end{align}
where
\begin{align}
g_n(\l)=\mathcal{L}_\xi^{-1}\left[\frac{(\xi^2+1)(2\xi)^n}{(1-\xi^2)^{n+1}}\right],
\end{align}
which can be explicitly computed for any $n.$

We can also obtain an explicit solution via a double integral as
follows. Making the substitution $s= \sinh( \ln\xi) = \frac12\left(
\xi - \frac1\xi \right),\  \xi > 1 $ in the equation

\begin{align}\label{neqn3}
\Phi(\sqrt{s^2+1}-s)=\sqrt{s^2+1}F(s).
\end{align}
we obtain
\begin{align}\label{neqn4}
 {1\over \xi} \Phi  \left({1\over \xi}\right) = \left( {1\over \xi^2} + 1\right) \int_0^\infty e^{-y \left( \xi - {1\over \xi} \right) }
f(2y) dy,\ \xi > 1.
\end{align}

To proceed, we use two operational relations for the Laplace
transform, which may be found in \cite{RK66}. Specifically, if $g
(s) =(\mathcal{L}h)(s)$, then  via the uniqueness property for the
Laplace transform
\begin{align}\label{neqn5}
{1\over s } g\left( {1\over s}\right) = \mathcal{L}\left[
\int_0^\infty J_0\left( 2\sqrt {\lambda y} \right) h(\lambda)
d\lambda \right] (s),
\end{align}
which is (21) on page 171 of \cite{RK66}. Further
\begin{align}\label{neqn6}
&\left(  {1\over s^2 }+ 1\right)  g\left(s-  {1\over s}\right) =
\mathcal{L} \left[   \int_0^y (y- \lambda)^{1/2} I_1\left( 2\sqrt
{\lambda (y - \lambda)} \right) h(\lambda) {d\lambda \over \sqrt
\lambda } \right.\nonumber\\&\left. -  \int_0^y (y- \lambda)^{- 1/2}
I_{-1}\left( 2\sqrt {\lambda ( y - \lambda)} \right) h(\lambda)
 \sqrt \lambda \  d\lambda \right] (s),
\end{align}
from formula (27), again p171 of \cite{RK66}. Hence, returning to
\eqnref{neqn4}, we cancel the Laplace transform via the uniqueness
property and take into account that $I_1(z)=I_{-1}(z)$ to obtain for
$ y >0$

\begin{align}\label{neqn7}
 \int_0^\infty J_0\left( 2\sqrt {\lambda y} \right) \varphi (\lambda) d\lambda
=  \int_0^y (y-2\lambda) \frac {I_1\left( 2\sqrt {\lambda (y -
\lambda)} \right)}{\sqrt{\lambda(y-\lambda)}}  f (2\lambda)
d\lambda.
\end{align}

However the left-hand side of \eqnref{neqn7} represents a variant of
the Hankel transform of  index zero, which admits the symmetric
inversion formula for $\varphi$ whose Mellin transform $\varphi^*(s)
\in L_1(\sigma),\ \sigma = \{ s \in \mathbb{C}: {\rm Re}( s )=
1/2\}$. (See the details in \cite{yal}, Ch.2, Section 2.2, Example
2.2). We remark here that the Mellin transform is defined by the
integral
\begin{align}\label{eqn5}
 f^*(s)= \int_0^\infty f(x) x^{s-1} dx.
\end{align}
Its inverse transform under certain conditions is given by the
integral
\begin{align}\label{eqn6}
f(x)= {1\over 2\pi i}  \int_{\gamma - i\infty}^{\gamma +i\infty}
f^*(s)  x^{-s} ds
\end{align}
See \cite{tit}, \cite{yal} for details of the Mellin transform.

Thus, inverting the Hankel transform, we arrive at the unique and
explicit solution of the integral equation \eqnref{eqn3one}, namely,
\begin{align}\label{neqn8}
\varphi(\lambda)=   \int_0^\infty  \int_0^{2y}  J_0\left( 2\sqrt
{\lambda y} \right)  \  \frac {I_1\left( \sqrt {u (2y - u)}
\right)}{\sqrt{u(2y- u)}}   (y-u) f (u) \  du dy,
\end{align}
$ \lambda >0.$   We note, that the interchange of the order of
integration in \eqnref{neqn8} is impossible because the
corresponding integral with respect to $y$ will be divergent.

Now we observe that the inequality
\begin{align}
1 < \Gamma(\nu+1)\left(\frac{2}{x}\right)^\nu I_\nu(x) < \cosh x<e^x
\end{align}
(see \cite{Luke72}) leads to
\begin{align}
\left|\int_0^{2y} \frac {I_1\left( \sqrt {u (2y - u)}
\right)}{\sqrt{u(2y- u)}}   (y-u) f (u) \  du\right|\leq
\frac{y}{2}e^y\int_0^{2y}|f(u)|du.
\end{align}
Zemanian constructed a Frechet space $\mathfrak{H}_\mu(\Real_+)$
characterized by three properties: Every $\phi\in
\mathfrak{H}_\mu(\Real_+)$ is a rapidly decreasing smooth function
on $\Real_+$, with an expansion of the form
\begin{align}
\phi(x)&=x^{\mu+1/2}\left(a_0+a_2x^2+\cdots+a_{2k}x^{2k}+R_{2k}(x)\right),
\end{align}
where
\begin{align}
a_{2k}=\frac{1}{k!2^k}\lim_{x\to0}(x^{-1}D)^k(x^{-\mu-1/2}\phi(x)),
\end{align}
and the remainder satisfies $(x^{-1}D)^kR_{2k}(x)=o(1)$ as $x\to
0^+.$ (See Lemma 5.2.1 of \cite{Zem68}). Zemanian proves that the
Hankel transform of index $\mu$ is an automorphism on
$\mathfrak{H}_\mu(\Real_+)$, \cite{Zem68}. We see then that if
$ye^y\int_0^{2y}|f(u)|du\in\mathfrak{H}_0(\Real_+)$, then
$\varphi\in \mathfrak{H}_0(\Real_+).$

We summarize our results by the following result.

\begin{thm}\label{prop1} Let $f $ be a given continuous function on
$\mathbb{R}_+$ such that $f(y)= O(y^{-1/2}),\ y \to \infty$. Let
$\varphi \in L_1(\mathbb{R}_+)$ and its Mellin transform
$\varphi^*(s) \in L_1(\sigma)$. Then $\varphi$ given by formula
\eqnref{neqn8} represents the unique solution of the integral
equation \eqnref{eqn3one}. Moreover if
$ye^yF(y)\in\mathfrak{H}_0(\Real_+)$, where
$F(y)=\int_0^{2y}|f(u)|du,$ and
\begin{align}
\tilde{f}(\l)= \int_0^\infty  \int_0^{2y}  J_0\left( 2\sqrt {\lambda
y} \right)\frac {I_1\left( \sqrt {z(2y - z)} \right)}{\sqrt{z(2y-
z)}}   (y-z) f (z)dz dy,
\end{align}
then
\begin{align}
u(y,t)= \frac12\int_0^\infty \tilde{f}(\l)\left[J_0(\sqrt{y^2+2\l
ye^t})+J_0(\sqrt{y^2+2\l ye^{-t}})\right] d\l,
\end{align}
 is a solution of
\begin{align}
u_{tt}=y^2u_{yy}+yu_y+y^2u,
\end{align}
satisfying the initial conditions $u(y,0)=f(y),
u_t(y,0)=0.$
\end{thm}

To obtain a solution satisfying $u(y,0)=0$ and $u_t(y,0)=g(y)$ let
\begin{align}\label{gsoln}
w(y,t)= \frac12\int_0^\infty \psi(\l)\left[J_0(\sqrt{y^2+2\l
ye^t})-J_0(\sqrt{y^2+2\l ye^{-t}})\right] d\l.
\end{align}
Then $w(y,0)=0$ and
\begin{align}
w_t(y,0)=-\int_0^\infty\psi(\l)\frac{\l y J_1\left(\sqrt{y^2+2 \l y
}\right)}{\sqrt{y^2+2\l y}}d\l.
\end{align}
We therefore seek to solve
\begin{align}\label{inteqn2}
-\int_0^\infty\l\psi(\l)\frac{ J_1\left(\sqrt{y^2+2 \l y
}\right)}{\sqrt{y^2+2\l y}}d\l=\frac{1}{y}g(y)=\tilde{g}(y).
\end{align}

This can actually be reduced to the case covered by Theorem
\ref{prop1}. In fact, assuming that $\lambda \psi(\lambda)$ is
absolutely continuous on $\mathbb{R}_+$ and its derivative
$\varphi(\lambda)= [\lambda \psi(\lambda) ]^\prime$ satisfies
conditions of Theorem \ref{prop1},  we integrate by parts in
\eqnref{inteqn2}, employing  the equality
$${y \over \sqrt{y^2+ 2\lambda y}} J_1(  \sqrt{y^2+ 2\lambda y} )= - {d\over d\lambda } J_0(\sqrt{y^2+ 2\lambda y} )$$
and obtain  the integral equation
 $$\int_0^\infty  [\lambda \psi(\lambda) ]^\prime J_0(\sqrt{y^2+ 2\lambda y} )  d\lambda = - g(y). $$
Hence our previous result gives
\begin{align}\label{neqn17}
\psi (\lambda)=  {1\over \lambda } \int_0^\lambda \int_0^\infty
\int_0^{2y}  J_0\left( 2\sqrt { y v} \right)  \ \frac {I_1\left(
\sqrt {u (2y - u)} \right)}{\sqrt{u(2y- u)}}   (u- y) g (u) \  du dy
dv,\
\end{align}
for $\lambda >0.$

\begin{thm}\label{prop2}
Suppose that $ye^yG(y)\in \mathfrak{H}_0(\Real_+)$,
$G(y)=\int_0^{2y}|g(z)|dz$. Then the equation
$u_{tt}=y^2u_{yy}+yu_y+y^2u$ with $u(y,0)=0$ and $u_t(y,0)=g(y)$ has
a solution
\begin{align}
u(y,t)=\frac12\int_0^\infty&\tilde{g}(\l)\left[J_0(\sqrt{y^2+2\l
ye^t})-J_0(\sqrt{y^2+2\l ye^{-t}})\right] d\l,
\end{align}
where
\begin{align}
\tilde{g}(\l)={1\over \lambda } \int_0^\lambda \int_0^\infty
\int_0^{2y}  J_0\left( 2\sqrt { y v} \right)  \ \frac {I_1\left(
\sqrt {z (2y - z)} \right)}{\sqrt{z(2y- z)}}   (z- y) g (z)dzdydv.
\end{align}

\end{thm}
Hence, provided that $f$ and $g$ satisfy the conditions of Theorems
\ref{prop1} and \ref{prop2}, a solution of our original problem may
be written
\begin{align}\label{gensoln}
u(y,t)&=\frac12\left( \int_0^\infty
\tilde{f}(\l)\left[J_0(\sqrt{y^2+2\l ye^t})+J_0(\sqrt{y^2+2\l
ye^{-t}})\right]d\l\nonumber \right.\\ & \left. + \int_0^\infty
\tilde{g}(\l)\left[J_0(\sqrt{y^2+2\l ye^t})-J_0(\sqrt{y^2+2\l
ye^{-t}})\right] d\l\right).
\end{align}
\subsubsection{The Direct Laplace Transform Approach}
It is worth showing how the Laplace transform method can be applied
to \eqnref{inteqn2}. This approach is computationally useful. On
page 60 of \cite{RK66} we find
\begin{align}
\int_0^\infty \frac{J_1(\sqrt{y^2+2\l y})}{\sqrt{y^2+2\l
y}}e^{-sy}dy=\frac{1-e^{-\l(\sqrt{s^2+1}-s)}}{\l}.
\end{align}
If the Laplace transform of $g(y)/y$ is $G$, then taking Laplace
transform of both sides of \eqnref{inteqn2}, we obtain
\begin{align}
-\int_0^\infty \psi(\l)(1-e^{-\l(\sqrt{s^2+1}-s)})d\l=G(s).
\end{align}
Which is the same as
\begin{align}
\Psi(\sqrt{s^2+1}-s)-\Psi(0)=G(s).
\end{align}
We have denoted the Laplace transform of $\psi$ by $\Psi$. Putting
$s=\sinh k$ and $k=-\ln\xi$ gives
\begin{align}
\Psi(\xi)-\Psi(0)=G\left(\frac12\left(\frac1\xi-\xi\right)\right).
\end{align}

So that
\begin{align}
\psi(\l)=\Psi(0)\delta(\l)+\mathcal{L}_\xi^{-1}\left[G\left(\frac12\left(\frac1\xi-\xi\right)\right)\right](\l).
\end{align}
The term $\Psi(0)\delta(\l)$ is in principle arbitrary, but making a
choice for $\psi$ will fix it, and in any case, it actually makes no
contribution to the solution since the two terms involving Bessel
functions cancel when evaluated at $\l=0.$

Once more, if
\begin{align}
g(y)=\sum_{n=0}^\infty \frac{1}{n!}g^{(n)}(0)y^n
\end{align}
then
\begin{align}
G\left(\frac12\left(\frac1\xi-\xi\right)\right)=\sum_{n=0}^\infty
g^{(n)}(0)\frac{(2\xi)^{n+1}}{(1-\xi^2)^{n+1}}
\end{align}
and so
\begin{align}
\psi(\l)=\Psi(0)\delta(\l)+\sum_{n=0}^\infty g^{(n)}(0)k_n(\l)
\end{align}
where
\begin{align}
k_n(\l)=\mathcal{L}_\xi^{-1}\left[\frac{(2\xi)^{n+1}}{(1-\xi^2)^{n+1}}\right].
\end{align}

Since the delta function terms make no contribution, the solution
can be written
\begin{align}
w(y,t)=\frac12\int_0^\infty \widehat{g}(\l)\left[J_0(\sqrt{y^2+2\l
ye^t})-J_0(\sqrt{y^2+2\l ye^{-t}})\right] d\l,
\end{align}
with $\widehat{g}(\l)=\sum_{n=0}^\infty g^{(n)}(0)k_n(\l).$

Finding examples using the direct Laplace transform approach is not
difficult.

\begin{example}
We consider the solution in the case when $f(y)=0$,
$g(y)=\frac14y^2{}_1F_2(\frac12;\frac32,2;-\frac{y^2}{4})-\frac12y.$
We find that
\begin{align}
\int_0^\infty
\frac{g(y)}{y}e^{-sy}dy=\frac{-1}{1+s+\sqrt{s^2+1}}=\Psi(\sqrt{s^2+1}-s)-\Psi(0).
\end{align}
Putting $s=\sinh k$ gives
\begin{align}
\frac{-1}{1+e^k}=\Psi(e^{-k})-\Psi(0)
\end{align}
and $\xi=e^{-k}$ leads to
\begin{align}
\frac{-\xi}{1+\xi}=\frac{1}{1+\xi}-1=\Psi(\xi)-\Psi(0),
\end{align}
and so
\begin{align}
\Psi(\xi)=\frac{1}{1+\xi}+\Psi(0)-1.
\end{align}

The inverse Laplace transform of both sides gives
\begin{align}
\psi(\l)=e^{-\l}+(\Psi(0)-1)\delta(\l).
\end{align}
The natural choice here is to take $\psi(\l)=e^{-\l}$ and then
$\Psi(0)-1=0.$ The solution our method gives for the initial value
problem is then
\begin{align}\label{e1}
u(y,t)=\frac12\int_0^\infty e^{-\l}\left[J_0(\sqrt{y^2+2\l
ye^t})-J_0(\sqrt{y^2+2\l ye^{-t}})\right] d\l.
\end{align}
This integral does not seem to be previously known. Although it is
easy to evaluate numerically, it is worthwhile to give an
evaluation.

In order to calculate \eqnref{e1} explicitly, we substitute the
series expansion for the Bessel function (see 9.1.10 of \cite{AS72}
) inside the integral and change the order of integration and
summation via the absolute and uniform convergence. Thus we obtain
\begin{align}\label{e2}
u_1(y,t)&= \int_0^\infty e^{-\lambda}  J_0(\sqrt{y^2+ 2\lambda y
e^t})  d\l\nonumber \\&= \sum_{k=0}^\infty { (-1)^k\over 4^k (k!)^2}
\int_0^\infty  e^{-\lambda} \left( y^2+ 2\lambda y e^t \right)^k
d\lambda .
\end{align}

Hence, the latter series becomes
\begin{align}\label{e3}
 \sum_{k=0}^\infty  { (-1)^k\over
4^k (k!)^2} \sum_{m=0}^k \  m!   \binom {k} {m}  y^ {2k-m} (2 e^t)^m
= \sum_{k=0}^\infty  { (-1)^k y^{2k} \over 4^k k! } \sum_{m=0}^k \ {
(2 e^t)^m \over (k-m)!  y^m}.
\end{align}
Making a straightforward substitution in the inner sum and then
changing the order of summation, we get
$$\sum_{k=0}^\infty  { (-1)^k y^{2k} \over 4^k k!  } \sum_{m=0}^k \  { (2 e^t)^m \over (k-m)!  y^m}=  \sum_{m=0}^\infty \sum_{k=m}^\infty  \left(-  {ye^t\over 2}\right)^{k-m}   \left(-  {y^2\over 4}\right)^m  {1 \over m!  k!  } $$
$$= \sum_{m=0}^\infty \sum_{k=0}^\infty  \left(-  {ye^t\over 2}\right)^k   \left(-  {y^2\over 4}\right)^m  {1 \over m!  (k+m)!  } . $$
The latter double series can be rewritten, recalling  the definition
of the Bessel function.  Hence the value of the integral $u_1(y,t)$
is
\begin{align}\label{e5}
u_1(y,t)= \sum_{k=0}^\infty  ( - 1 )^k e^{kt} J_k(y).
\end{align}

The series \eqnref{e5} has  no simple expression, but can be written
in terms of the Lommel function of two variables \cite{prud},
\begin{align}\label{e6}
U_\nu( x,y)=       \sum_{k=0}^\infty  ( - 1 )^k \left({x\over
y}\right)^{2k+\nu}   J_{2k+\nu} (y),\  x,y > 0.
\end{align}

Thus appealing to relation (5.7.5.1) in \cite{prud}, we write
\eqnref{e5} in the form
$$u_1(y,t)= U_0( - e^t i y, y) + i U_1( - e^t i y, y),$$
where $i=\sqrt{-1}$ and, correspondingly,  the value of the integral
\eqnref{e1} as
\begin{align}
 u(y,t)&= \frac{\left[ U_0( - e^t i y, y) - U_0( - e^{-t} i y,
y)\right]  +
  \left[ U_1( - e^t i y, y)-  U_1( - e^{-t} i y, y)\right]}{2}.
\end{align}
\end{example}

The interested reader can produce many examples. For instance, the
choice $f(y)=0, g(y)=yJ_0(y)$ leads to $\psi(\l)=2\cos(\l)$, for
which $\Psi(0)=0$.

\subsubsection{The Case $A\neq 0$ and $B=0$}.

Now let us briefly consider the equation
\begin{align}
u_{tt}=u_{xx}+\frac{1}{A^2e^{2x}}u,\ x\in \mathbb{R}.
\end{align}
Without loss of generality we take $A=1$. Now set $y=e^{-x}.$ This
again produces the equation
\begin{align}
u_{tt}=y^2u_{yy}+yu_y+y^2u,\ y>0.
\end{align}
The analysis then proceeds exactly as in the previous case.

\section{The Case When $A$ and $B$ are both nonzero}
Now consider the equation $u''(x)+\frac{1}{(Ae^{\omega
x}+Be^{-\omega x})^2}u=0. $ Suppose that $\omega, AB>0$ and let
$\alpha=\frac{\sqrt{AB(1+AB)}}{AB}.$ Then with
$\beta=\frac{1-\alpha}{2},$
\begin{align*}
u(x)=\left(1+\frac{A}{B}e^{2\omega
x}\right)^{\beta}{}_2F_1\left(\beta,\beta;2\beta;-\left(1+\frac{A}{B}e^{2\omega
x}\right)\right),
\end{align*}
is a stationary solution of the PDE
$u_{tt}=u_{xx}+\frac{1}{(Ae^{\omega x}+Be^{-\omega x})^2}u$. Here
${}_2F_1$ is Gauss' hypergeometric function, (Chapter 15 of
\cite{AS72}). Applying the symmetry we obtain the new, non
stationary solution

\begin{align*}
K(x,t;\ep)&=\left(\frac{1+\frac{A}{B}e^{2\omega x}}{1+2A\ep\omega
e^{\omega(x+t)}}\right)^{\beta}{}_2F_1\left(\beta,\beta;2\beta;\frac{-(1+\frac{A}{B}e^{2\omega
x})}{1+2A\ep\omega e^{\omega(x+t)}}\right).
\end{align*}
We assume that $\beta$ is not a negative integer here. We can
immediately construct an integral operator mapping functions to
solutions by setting

\begin{align}
u(x,t)&=\frac12 \int_0^\infty
\varphi(\ep)\left[K(x,t;\ep)+K(x,-t;\ep)\right]d\ep
\end{align}
for suitable $\varphi.$ This will satisfy $u_t(y,0)=0$ and
\begin{align}\label{generalinteqn}
u(x,0)=\int_0^\infty \varphi(\ep)\left(\frac{1+\frac{A}{B}e^{2\omega
x}}{1+2A\ep\omega e^{\omega
x}}\right)^{\beta}{}_2F_1\left(\beta,\beta;2\beta;\frac{-(1+\frac{A}{B}e^{2\omega
x})}{1+2A\ep\omega e^{\omega x}}\right)d\ep.
\end{align}

If we take
\begin{align}
u(x,t)=\frac12\int_0^\infty \psi(\ep)[K(x,t;\ep)-K(x,-t;\ep)]d\ep,
\end{align}
we obtain a solution satisfying $u(x,0)=0$ and
\begin{align}
u_t(x,0)=\int_0^\infty \psi(\ep)\frac{\partial}{\partial
t}K(x,t;\ep)|_{t=0}d\ep.
\end{align}
If we require $u_t(x,0)=g(x)$, we have an integral equation for
$\psi$.

We will not attempt to solve these difficult integral equations
here. Rather we look at the special case  $\omega=1$ and $A=B=1/2.$
Then we obtain the solution
\begin{align*}
u(x,t)&=c_1P_{\nu}\left(\frac{2
   \epsilon  e^{t+x}+e^{2 x}-1}{e^{2 x}+1}\right)+c_2Q_{\nu}\left(\frac{2
   \epsilon  e^{t+x}+e^{2 x}-1}{e^{2 x}+1}\right),
\end{align*}
for the equation
\begin{align}\label{hypex2}
u_{tt}&=u_{xx}+(\text{sech}^2x) u.
\end{align}
In this $ \nu=\frac{1}{2} \left(\sqrt{5}-1\right). $ Suppose we set
$c_1=0,c_2=1.$ Then we can form a solution of the PDE by setting
\begin{align}\label{KLinc}
u(x,t)=\int_0^\infty \varphi(\ep)Q_{\nu}\left(\frac{2
   \epsilon  e^{t+x}+e^{2 x}-1}{e^{2 x}+1}\right)d\ep.
\end{align}
Different domains of integration can also be considered. The
operator defined by \eqnref{KLinc} maps functions $\varphi$ to
solutions of \eqnref{hypex2}. We will show how this operator can be
used to solve  the initial value problem
\begin{align}
u_{tt}&=u_{xx}+(\text{sech}^2x) u\\
u(x,0)&=f(x),\\
u_t(x,0)&=0.
\end{align}

We suppose that $f$ lies in some suitable function space, which will
be discussed below. Now we remark that standard Laplace transform
results and formula 23(3) on page 337 of \cite{RK66} gives the
result
\begin{align*}
\int_0^\infty \frac{\sqrt{\frac{\pi }{2}} e^{-e^{-t} \xi  \sinh
    x -\frac{t}{2}} I_{\nu +\frac{1}{2}}\left(e^{-t} \xi
   \cosh  x \right)}{\sqrt{\xi  \text{sech} x }}e^{-\ep\xi}d\xi=Q_{\nu}\left(\frac{2
   \epsilon  e^{t+x}+e^{2 x}-1}{e^{2 x}+1}\right).
\end{align*}
Thus formula \eqnref{KLinc} can be rewritten
\begin{align}\label{solnintform}
w(x,t)&=\int_0^\infty\int_0^\infty\varphi(\ep)\frac{\sqrt{\frac{\pi
}{2}} e^{-e^{-t} \xi  \sinh
    x -\frac{t}{2}} I_{\nu +\frac{1}{2}}\left(e^{-t} \xi
   \cosh  x \right)}{\sqrt{\xi  \text{sech} x }}e^{-\ep\xi}d\xi
   d\ep\nonumber \\
   &=\int_0^\infty\Phi(\xi)\frac{\sqrt{\frac{\pi
}{2}} e^{-e^{-t} \xi  \sinh
    x -\frac{t}{2}} I_{\nu +\frac{1}{2}}\left(e^{-t} \xi
   \cosh  x \right)}{\sqrt{\xi  \text{sech} x }} d\xi,
\end{align}
in which $\Phi$ is the Laplace transform of $\varphi.$ We obviously
require suitable decay for $\varphi$ in order to justify the use of
Fubini's Theorem. Thus from the symmetry solution we have obtained
an operator which maps functions $\Phi$ from a suitable test space
to solutions of \eqnref{hypex2}.

If we let $\xi\to i\xi$ and use linearity we obtain the solutions
\begin{align}
w_1(x,t)&=\frac{ \cos(e^{-t} \xi  \sinh
    x)e^{ -\frac{t}{2}} J_{\nu +\frac{1}{2}}\left(e^{-t} \xi
   \cosh  x \right)}{\sqrt{\xi  \text{sech} x }}\\
w_2(x,t)&=\frac{ \sin(e^{-t} \xi  \sinh
    x)e^{ -\frac{t}{2}} J_{\nu +\frac{1}{2}}\left(e^{-t} \xi
   \cosh  x \right)}{\sqrt{\xi  \text{sech} x }}.
\end{align}
From these we obtain the solution of \eqnref{hypex2} given by
\begin{align*}
u(x,t)&=\sqrt{\pi\over 2}\int_0^\infty \Phi(\xi)\frac{ \cos(e^{-t}
\xi  \sinh
    x)e^{ -\frac{t}{2}} J_{\nu +\frac{1}{2}}\left(e^{-t} \xi
   \cosh  x \right)}{\sqrt{\xi  \text{sech} x }}d\xi \\
   &+\sqrt{\pi\over 2}\int_0^\infty \Psi(\xi)\frac{ \sin(e^{-t} \xi  \sinh
    x)e^{ -\frac{t}{2}} J_{\nu +\frac{1}{2}}\left(e^{-t} \xi
   \cosh  x \right)}{\sqrt{\xi  \text{sech} x }}d\xi.
\end{align*}

Notice that time reversal is also a symmetry and so
\begin{align*}
\bar{u}(x,t)&=\sqrt{\pi\over 2}\int_0^\infty \Phi(\xi)\frac{
\cos(e^{t} \xi  \sinh
    x)e^{ \frac{t}{2}} J_{\nu +\frac{1}{2}}\left(e^{t} \xi
   \cosh  x \right)}{\sqrt{\xi  \text{sech} x }}d\xi \\
   &+\sqrt{\pi\over 2}\int_0^\infty \Psi(\xi)\frac{ \sin(e^{t} \xi  \sinh
    x)e^{\frac{t}{2}} J_{\nu +\frac{1}{2}}\left(e^{t} \xi
   \cosh  x \right)}{\sqrt{\xi  \text{sech} x }}d\xi,
\end{align*}
is again a solution of the PDE. We will therefore let

\begin{align*}
u(x,t)&=\frac12\sqrt{\pi\over 2}\int_0^\infty \Phi(\xi)\left[\frac{
\cos(e^{t} \xi \sinh
    x)e^{ \frac{t}{2}} J_{\nu +\frac{1}{2}}\left(e^{t} \xi
   \cosh  x \right)}{\sqrt{\xi  \text{sech} x }}\right.\\&\left.\quad +\frac{ \cos(e^{-t} \xi  \sinh
    x)e^{ -\frac{t}{2}} J_{\nu +\frac{1}{2}}\left(e^{-t} \xi
   \cosh  x \right)}{\sqrt{\xi  \text{sech} x }}\right]d\xi.
\end{align*}

It is straightforward to see that

\begin{align}
u_{tt}=u_{xx}+(\text{sech}^2x) u,
\end{align}
\begin{align}
u(x,0)=\sqrt{\pi\over 2} \int_0^\infty \Phi(\xi)\frac{ \cos(  \xi
\sinh
    x) J_{\nu +\frac{1}{2}}\left( \xi
   \cosh  x \right)}{\sqrt{\xi  \text{sech} x }}d\xi
\end{align}

and
\begin{align}
u_t(x,0)=0.
\end{align}

If we absorb the $\sqrt{\xi}$ into $\Phi$, which is after all
arbitrary, we see that if we can solve the integral equation
\begin{align}
f(x)\sqrt{\text{sech} x}=\sqrt{\pi\over 2}\int_0^\infty  \Phi(\xi)
\cos( \xi \sinh
    x) J_{\mu}\left( \xi
   \cosh  x \right)  d\xi,
\end{align}
for $\Phi$, with $\mu=\nu +\frac{1}{2}$, then we have a solution of
our Cauchy problem. This is not a standard integral equation and
does not seem to appear anywhere in the existing literature. The
solution is more involved than the previous integral equations we
solved. Consequently, we will discuss its solution in the following
section.

\section{Solving the Integral Equation}

Let $x \in \mathbb{R}_+$ and consider the following non-convolution
integral equations of the first kind
\begin{align}\label{eqn1}
\sqrt{\pi\over 2} \int_0^\infty \varphi(\xi) \cos(\xi\sinh x)
J_\nu(\xi\cosh x) d\xi &= f(x),
\end{align}
\begin{align}\label{eqn2}
 \sqrt{\pi\over 2} \int_0^\infty \varphi(\xi) \sin (\xi\sinh x)
J_\nu(\xi\cosh x) d\xi &= f(x),
\end{align}
where $f(x) $ is a given function and $\varphi (\xi)$ is to be
determined and $J_\nu(z),\  \nu \in \mathbb{C}$ is the Bessel
function of the first kind \cite{AS72}. Though only the solution of
the first equation is needed for our analysis, we present the
solution of both.

The key ingredients to solve integral equations \eqnref{eqn1} and
\eqnref{eqn2} are the following integrals (see \cite{GR2000},
6.699(2) and 6.699(1), p.723, with $\l=s-1$),
\begin{align}\label{eqn3}
\int_0^\infty x^{s-1} \cos(bx) J_\nu(cx) dx&= \frac{}{}
\frac{2^{s-1}\Gamma\left(\frac{s+\nu}{2}\right)}{c^s\Gamma(1+
\frac{(\nu-s)}{2})}\ {}_2F_1 \left( {s+\nu\over 2},  {s- \nu\over
2}; {1\over 2}; {b^2\over c^2}\right),
\end{align}
where $c >b,\  -  {\rm Re}( \nu) < {\rm Re }( s) < 3/2, \ \Gamma(z)$
is the Euler gamma function, see Chapter 6 of \cite{AS72}. We also
have
\begin{align}\label{eqn4}
 \int_0^\infty x^{s-1} \sin (bx) J_\nu(cx) dx&=    \frac{2^{s} b \Gamma\left(\frac{s+\nu+1}{2}\right)}{c^{s+1}\Gamma\left(\frac{\nu-s+1}{2}\right)}
\nonumber \\&\times  {}_2F_1 \left( {s+\nu+1\over 2},  {s+1 -
\nu\over 2}; {3\over 2}; {b^2\over c^2}\right).
\end{align}
As we see the left-hand sides of \eqnref{eqn3} and \eqnref{eqn4}
represent the Mellin transforms of the kernels in equations
\eqnref{eqn1} and \eqnref{eqn2}. Recall that for two functions $f,
g$ the generalized Mellin-Parseval equality holds (see  \cite{tit},
\cite{yal} for the details)
\begin{align}\label{eqn7}
\int_0^\infty f(xy) g(y) dy = \frac{1}{ 2\pi i} \int_{\gamma-
i\infty}^{\gamma+i\infty} f^*(s) g^*(1-s) x^{-s} ds,\ x>0.
\end{align}
The Fourier cosine and sine transforms of the integrable function
$f$ are  defined by  the formulas
\begin{align}\label{eqn8}
(F_{c\brace s} f)( x) = \sqrt{\frac{2}{\pi}} \int_0^\infty { \cos
(xy) \brace \sin (x y) } f(y)dy,
\end{align}
see \cite{tit}.

This equation does not have a unique solution. However if we impose
certain extra conditions on $\varphi$, uniqueness can be obtained.
We present the general case first.

\begin{thm}\label{Main}  Let ${\rm Re}(\nu) > -1/2, \varphi(\xi) \in  S
(\mathbb{R}_+),\ e^{3x/2} \  f (x)  \in L_1(\mathbb{R}_+), \ $ and
$y^{-\nu-1/2} (F_c g)( y) \in S (\mathbb{R}_+)$,\ where $ g(x)=
(x^2+1)^{-\nu/2}  \  f \left(\sinh^{-1}  x\right)$. Let also the
Mellin transform $\varphi^*$ of $\varphi$ satisfy the condition

$$\varphi^*(s)  |s|^{-\nu}  e^{(\pi/2 -\delta) |s|} \in L_1(\sigma),\  \delta \in [0, \pi/2) ,$$
where  $\sigma = \left\{ s\in \mathbb{C}, \  s = {1\over 2}+i\tau,\
\tau \in \mathbb{R}\right\}$ and
\begin{align}\label{eqn9}
h_{c}(\tau)&=   \varphi^*(1/2 +i\tau)  \Gamma(1/2- i\tau+\nu)  +
\varphi^*(1/2 -i\tau) \Gamma(1/2+ i\tau+\nu) \\&=
O(\tau^2),\nonumber
\end{align}
$\tau \to 0.$ Then  the integral equation \eqnref{eqn1} has the
following solutions \
\begin{align}\label{eqn10}
&\varphi(\xi) = - {1\over 4\pi i} \int_\sigma {\rho(s)\over
\Gamma(1-s+\nu) }\ \xi^{-s} ds+ { 2\xi^{-(1+\nu)} \over \pi\sqrt \pi } \times \nonumber \\
& \int_{-\infty}^\infty e^{(\nu+1/2)u} \left[  \left({1\over \sqrt
2}+\sqrt 2\ \nu\right) \sin\left({\pi\over 2} \left(\nu+{1\over 2}
\right) - e^{u}
\xi^{-1}\right) +  \cos (\pi\nu/2)\right.\nonumber\\
& \times \left[  {  e^{u} \over \sqrt 2 \xi}\ \left( \cos\left(e^{u}
\xi^{-1} \right)  -  {1\over \sqrt 2} \sin \left(e^{u} \xi^{-1}
\right) \right)  -  \sin \left(e^{u} \xi^{-1} \right)  \left(1+
{1\over  \sqrt 2}\right) \right]\nonumber \\
&\left.+  \sin \left({\pi\nu \over 2} \right)  \sin \left(e^{u}
\xi^{-1} \right) \left( 1- {e^{u} \over  2 \xi}\right)\right]\nonumber\\
&\times  \int_0^\infty \   y^{-1/2-\nu} \cos (y\sinh u)
\int_0^\infty \cos( y\sinh x)  (\cosh x)^{1-\nu} f(x) \ dx dy du,
\end{align}
where $\xi>0$ and depending on an arbitrary function  $\rho(s)$,
which is odd on $\sigma$, i.e.
\begin{align}\label{eqn11}
\rho(s)= -\rho(1-s),  s= 1/2+i\tau.
\end{align}
\end{thm}

\begin{proof}  Let $\sigma= \left\{ s\in \mathbb{C}, \  s = {1\over 2}+i\tau,\ \tau \in \mathbb{R}\right\}$.
Since via integration by parts in the integral \eqnref{eqn5} it is
not difficult to verify that  the Mellin transform $\varphi^*(s)$ of
$\varphi \in S (\mathbb{R}_+)$  belongs to $L_1(\sigma),$  we have
by the inversion formula \eqnref{eqn6}
\begin{align}\label{eqn12}
\varphi(\xi)=  {1\over 2\pi i}  \int_{\sigma} \varphi^*(s)  \xi^{-s}
ds.
\end{align}

Substituting this integral into \eqnref{eqn1}, we can proceed if we
may change the order of integration.  But since ( see 9.1.7 of
\cite{AS72}) $J_{\nu}(x)= O(x^\nu),\ x \to 0$ and ${\rm Re}(\nu)
> -1/2$, this will be possible if we show (see, for instance, in
\cite{yal}, Section 2.1) that there exists a positive absolute
constant $C$ such that any fixed $x >0$ and for almost all $E_1, E_2
> 1$ and all $\tau  \in \mathbb{R}$
$$ \left| \int_{E_1}^{E_2}  \xi^{-1/2- i\tau} \cos(\xi\sinh x) J_\nu(\xi\cosh x) d\xi\right| < C.$$
The latter inequality follows, employing the asymptotic behavior of
the Bessel function at infinity (see 9.2.1 of \cite{AS72})
$$J_{\nu}(x)\sim   \sqrt{ {2\over \pi x}} \cos\left( x - {\pi\nu\over 2}\right), \ x \to \infty$$
and from  the uniform convergence with respect to $\tau$ of the
following integrals
\begin{align}
\int_{E_1}^{E_2}  \xi^{-1 - i\tau} \cos\left(\xi(\sinh x\pm \cosh x)
\mp {\pi\nu\over 2}\right) d\xi.
\end{align}
This fact can be verified with the use of the mean value theorems
and integration by parts,  and we omit the details.  Hence with the
use of the Boltz formula,
\begin{align}
{}_2F_1(a,b,c;z)=(1-z)^{-a}{}_2F_1\left(a,c-b;c;\frac{z}{1-z}\right)
\end{align}
(see 15.3.4 of \cite{AS72}), in the right-hand side of
\eqnref{eqn3}, the equation \eqnref{eqn1} takes the form
\begin{align}\label{eqn13}
\sqrt{\pi\over 2}  {\cosh^\nu x \over 2\pi i}&  \int_{\sigma} \frac{
\Gamma\left(\frac{1-s+\nu}{2}\right)}{2^{s}\Gamma\left(\frac{1+\nu+s}{2}\right)}
{}_2F_1 \left( {1-s+\nu\over 2},  {\nu + s\over 2}; {1\over 2};-
\sinh^2 x \right)   \varphi^*(s)  ds\nonumber \\& = f(x).
\end{align}

To proceed, we make the the substitution $s= 1/2+i\tau$ in the
integral \eqnref{eqn13} to obtain
\begin{align}\label{eqn13a}
\sqrt{\pi\over 2}  {\cosh^\nu x \over 2\pi  }& \int_{-\infty}^\infty
\frac{ \Gamma\left(\frac{1/2-i\tau+\nu}{2}\right){}_2F_1 \left(
{1/2-i\tau+\nu\over 2},  {1/2+\nu + i\tau\over 2}; {1\over 2};-
\sinh^2 x \right)
}{2^{1/2+i\tau}\Gamma\left(\frac{3/2+\nu+i\tau}{2}\right)}\nonumber
\\&\times \varphi^*(1/2+i\tau)  d\tau = f(x).
\end{align}
Next we split the integral as $\int_{-\infty}^\infty
=\int_{-\infty}^0+\int_{0}^\infty$ and note the symmetry
${}_2F_1(a,b;c;z)={}_2F_1(b,a;c;z)$. The substitution $\tau\to-\tau$
in the first integral and the duplication formula for the gamma
function,
$$\Gamma(2z)=\frac{1}{\sqrt{2\pi}}2^{2z-1/2}\Gamma(z)\Gamma(z+1/2),$$
(6.1.18 of \cite{AS72}) leads to the new form of \eqnref{eqn1}
\begin{align}\label{eqn14}
\int_{0}^\infty  h_{c}(\tau)\frac{ {}_2F_1 \left( {1/2+\nu+i\tau
\over 2}, {1/2+\nu -i\tau \over 2}; {1\over 2};  - \sinh^2 x
\right)}{ \Gamma\left({3/2+ \nu+i\tau\over 2}\right)
\Gamma\left({3/2+ \nu- i\tau\over 2}\right) }  d\tau = 2^{3/2+\nu}
\frac{f(x)}{\cosh^{\nu}( x)}.
\end{align}

Moreover, the integral \eqnref{eqn14} converges absolutely by virtue
of the Stirling asymptotic formula for the gamma function ( see
6.1.37 of \cite{AS72}) and the behavior at infinity in $\tau$ of the
hypergeometric function (cf. \cite{yak}, Theorem 1.12)

$${}_2F_1 \left( {1/2+\nu+i\tau \over 2},  {1/2+\nu -i\tau \over 2}; {1\over 2};  - \sinh^2 x \right) = O(1), \tau \to \infty.$$
In the mean time, the latter hypergeometric function has the
representation (cf. \cite{prud}, Vol. 2, relation (2.16.21.1),
\cite{yak}, formula (1.101))

\begin{align}\label{eqn15}
&{}_2F_1 \left( {1/2+\nu+i\tau \over 2}, {1/2+\nu -i\tau \over 2};
{1\over 2};  - \sinh^2 x \right)\nonumber =
\\&\times\frac{2^{3/2-\nu} }{\Gamma\left(\frac{ 1/2+ \nu+i\tau}{2}\right) \Gamma\left(\frac{1/2+
\nu- i\tau}{2}\right)} \int_0^\infty y^{\nu -1/2} \cos(y\sinh x)
K_{i\tau} (y) dy,
\end{align}
${\rm Re}( \nu) > -1/2,$ where $K_{i\tau}(x)$ is the modified Bessel
function of the third kind or the Macdonald function \cite{AS72},
and it is the kernel of the Kontorovich-Lebedev transform
\cite{yak}.

Substituting the right -hand side of \eqnref{eqn15} into the
left-hand side of \eqnref{eqn14},  we change the order of
integration by the Fubini theorem, which is justified via the
absolute convergent of the iterated integrals. Hence after the use
of the duplication formula for the gamma-function and the
substitution $u=\sinh x$, we arrive at the following integral
equation

\begin{align}\label{eqn16}
\int_0^\infty \int_{0}^\infty &{h_{c}(\tau) \ K_{i\tau} (y) y^{\nu
-1/2}  \cos(u y) \over
\left|\Gamma\left(\frac12+\nu+i\tau\right)\right|^2} d\tau dy  =
2\pi\frac{ f \left(\sinh^{-1} \ u\right)}{(u^2+1)^{\nu/2}},
\end{align}
since
$\Gamma(z)\Gamma(\bar{z})=\Gamma(z)\overline{\Gamma(z)}=|\Gamma(z)|^2.$

Further, by the definition of $h_c(\tau)$,  the Stirling asymptotic
formula for the Gamma function, the uniform inequality for the
Macdonald function \cite{yak}
$$|K_{i\tau} (y)| \le e^{-\delta \tau /2} K_0(y \cos\delta),\ \delta \in \left[0,  {\pi\over 2}\right),$$
and our assumption $\varphi^*(s)  |s|^{-\nu}  e^{(\pi/2 -\delta)
|s|} \in L_1(\sigma)$, one can apply the Fourier cosine transform to
both sides of the equality \eqnref{eqn16}. Hence,  returning to the
original variable,  we arrive at the equation
\begin{align}\label{eqn17}
&\int_{0}^\infty  {h_{c}(\tau)\ K_{i\tau} (y) \over
\left|\Gamma\left(\frac12+\nu+i\tau\right)\right|^2}  d\tau  &= 4
 \int_0^\infty \frac{ \cos( y\sinh x)}{y^{-1/2+\nu}}  (\cosh x)^{1-\nu} f(x)
dx.
\end{align}

The left-hand side of\eqnref{eqn17} is related to the
Kontorovich-Lebedev transform \cite{yak}. Recall that for suitable
functions $g$, we have
\begin{align}\label{KLident}
g(x)=\frac{2}{\pi^2x}\int_0^\infty \int_0^\infty \tau\sinh(\pi
\tau)K_{i\tau}(x)K_{i\tau}(y)g(y)dyd\tau.
\end{align}
If we denote the right hand side of \eqnref{eqn17} by $\phi(y)$ and
introduce
\begin{align}
H_c(\tau)=\frac{\pi^2h_c(\tau)}{2\tau\sinh(\pi\tau)
\left|\Gamma\left(\frac12+\nu+i\tau\right)\right|^2},
\end{align}
then
\begin{align}
\phi(y)=\frac{2}{\pi^2}\int_0^\infty \tau\sinh(\pi
\tau)H_c(\tau)K_{i\tau}(y)d\tau.
\end{align}
But the conditions of the theorem permit us to invert the
Kontorovich-Lebedev transform, so by \eqnref{KLident} we have
\begin{align}
H_c(\tau)=\int_0^\infty \frac{\phi(y)}{y}K_{i\tau}(y)dy.
\end{align}
Thus we obtain

\begin{align}\label{eqn18}
h_{c}(\tau) &=  {8\over \pi^2} \tau\sinh(\pi\tau)
\left|\Gamma\left(\frac12+\nu+i\tau\right)\right|^2\int_0^\infty
y^{-1/2-\nu} K_{i\tau} (y)\nonumber \\&\times   \int_0^\infty
\frac{\cos( y\sinh x) }{(\cosh x)^{ \nu-1}}f(x)dx dy.
\end{align}

In the right-hand side of \eqnref{eqn18}, we employ the Parseval
equality for the Fourier cosine transform \cite{tit}, which is
allowed via conditions of the theorem. Then using relation
(2.16.14.1) in \cite{prud}, Vol. 2 and simple substitutions, the
latter equality \eqnref{eqn18} becomes
\begin{align}\label{eqn19}
h_{c}(\tau) =  {8\sqrt 2 \over \pi\sqrt\pi}
\tau\sinh\left({\pi\tau\over 2}\right)
\left|\Gamma\left(\frac12+\nu+i\tau\right)\right|^2\int_0^\infty \cos (\tau u)\nonumber \\
\times  \int_0^\infty \frac{ \cos (y\sinh u)}{  y^{\nu-1/2}}
\int_0^\infty \frac{\cos( y\sinh x) }{(\cosh x)^{ \nu-1}}f(x) \ dx
dy du.
\end{align}

Now, recalling  the definition of $h_{c}(\tau)$ in terms of
$\varphi^*$, let us write the equation \eqnref{eqn19} with respect
to the variable $s= 1/2+i\tau$, dividing first both of its sides by
$\Gamma(1/2+\nu +i\tau)$ and appealing to the reflection formula for
the Gamma function.  After making the change of variable $t= e^u$,
we have
\begin{align}\label{eqn20}
 &\varphi^*(s) {\Gamma(1-s+\nu)\over \Gamma(s+\nu)}  + \varphi^*(1-s) =  {8\sqrt 2 \over \sqrt\pi}
\frac{\Gamma( 1-s+\nu)}{ \Gamma\left( \frac{ s-1/2}{2}\right)
\Gamma\left(\frac{1/2-s}{2}\right)}\times\nonumber
\\&  \int_1^\infty   \left[ t^{s-\frac32} + \frac{1}{t^{s+\frac12}} \right]
\int_0^\infty \frac{  \cos (\frac12y(t-\frac1t) )}{ y^{1/2+\nu}}
\int_0^\infty \frac{\cos( y\sinh x) }{(\cosh x)^{ \nu-1}} f(x)   dx
dy dt,\ \
\end{align}
where $s= 1/2+i\tau, {\rm Re}(\nu) > - 1/2.$ Further, denoting the
iterated  integral in the right-hand side of \eqnref{eqn20} by
\begin{align}\label{eqn21}
F(s) &=  {8\sqrt 2 \over \sqrt\pi}   \int_1^\infty  \left[
t^{s-\frac32} + \frac{1}{t^{s+\frac12}} \right]  \int_0^\infty
\frac{  \cos (\frac12y(t-\frac1t) )}{ y^{1/2+\nu}} \nonumber
\\&\times \int_0^\infty \frac{\cos( y\sinh x) }{(\cosh x)^{
\nu-1}}f(x)  dx dy dt,
\end{align}
equation \eqnref{eqn20} becomes
\begin{align}\label{eqn22}
 \varphi^*(s) {\Gamma(1-s+\nu)\over \Gamma(s+\nu)}  + \varphi^*(1-s) =
\frac{ F(s)\  \Gamma( 1-s+\nu)}{ \Gamma\left( \frac{
s-1/2}{2}\right) \Gamma\left(\frac{1/2-s}{2}\right)}.
\end{align}
Note that $F(s)=F(1-s).$

Hence,
\begin{align}
\varphi^*(1-s) {\Gamma(s+\nu)\over \Gamma(1-s+\nu)}  + \varphi^*(s)
= \frac{ F(1- s)\  \Gamma( s+\nu)}{ \Gamma\left( \frac{
s-1/2}{2}\right) \Gamma\left(\frac{1/2-s}{2}\right)},
\end{align}
or
\begin{align}\label{third}
\varphi^*(1-s)  + \varphi^*(s) {\Gamma(1-s+\nu)\over \Gamma(s+\nu)}
=    {F(1-s)\  \Gamma(1-s+\nu)\over \Gamma\left( \frac{
s-1/2}{2}\right) \Gamma\left(\frac{1/2-s}{2}\right)}.
\end{align}

Therefore, adding equations \eqnref{eqn22} and \eqnref{third}, we
find
\begin{align}\label{eqn532}
2 \varphi^*(1-s)  + 2 \varphi^*(s) {\Gamma(1-s+\nu)\over
\Gamma(s+\nu)}  =   {\left[ F(1-s)  + F(s)\right]\
\Gamma(1-s+\nu)\over  \Gamma\left( \frac{ s-1/2}{2}\right)
\Gamma\left(\frac{1/2-s}{2}\right)}.
\end{align}

Hence
$$\left[ 2 \varphi^*(1-s) - {  F(s) \  \Gamma(1-s+\nu)\over \Gamma\left(
\frac{ s-1/2}{2}\right) \Gamma\left(\frac{1/2-s}{2}\right)}\right]
\Gamma(s+\nu) = -   \Gamma(1-s+\nu) $$
$$\times \left[  2 \varphi^*(s) -  {F(1-s) \Gamma(s+\nu)\over \Gamma\left(
\frac{ s-1/2}{2}\right) \Gamma\left(\frac{1/2-s}{2}\right)} \right]
.$$

If we denote the right hand side of the latter algebraic equation by
$\rho(s)$, we see that it satisfies the parity property
\eqnref{eqn11}. Choice of either $\varphi^*(s)$ or $\varphi^*(1-s)$
fixes the other. Hence final solutions of the equation \eqnref{eqn1}
can be expressed in terms of the inverse Mellin transform
\eqnref{eqn6}

\begin{align}\label{eqn24}
 \varphi(\xi) &= - {1\over 4\pi i} \int_\sigma
{\rho(s)\xi^{-s}\over \Gamma(1-s+\nu) }\  ds  +  {1\over 4\pi i}
\int_\sigma { F(1-s) \  \Gamma\left(1-s+\nu\right)\over \Gamma\left(
\frac{s-1/2}{2}\right) \Gamma\left(\frac{1/2-s}{2}\right)}\xi^{-s}
ds,
\end{align}
where $\xi > 0$ and $\rho(s)$ is an arbitrary function, which is odd
on the line $\sigma$, i.e. satisfies \eqnref{eqn11}. The integrals
converge via the relation \eqnref{eqn20} and properties of
$\varphi^*(s)$ that we have assumed.   Our final goal is to
calculate the second integral in the right -hand side of
\eqnref{eqn24}. To do this, we will employ the generalized Parseval
equality of type \eqnref{eqn7} (see \cite{mar}, Chapter 7, Theorem
23), because the gamma-ratio behaves on $\sigma$ as
$${ \Gamma(1-s+\nu)\over \Gamma( (s-1/2)/2) \Gamma( (1/2-s)/2)}=  O(|s|^{1+ {\rm Re}( \nu)}), \   |{\rm Im} s| \to \infty,$$
and the inverse Mellin transform generally does not exist since  the
corresponding integral \eqnref{eqn6} may diverge. So, recalling
\eqnref{eqn21}, we observe that according to  conditions of the
theorem, the iterated integral
$$ \int_0^\infty \   y^{-1/2-\nu} \cos (y(t-1/t)/2)  \int_0^\infty  \cos( y\sinh x)  (\cosh x)^{1-\nu} f(x) \ dx dy$$
is a function of $(t-1/t)/2$ from the Schwartz space. Therefore, due
to the absolute and uniform convergence of the corresponding
integral \eqnref{eqn21},  $F(1-s)$  is analytic in the right half
-plane ${\rm Re( s)} > 1/2$.  Hence, returning to \eqnref{eqn24} and
using the duplication formula for the gamma function
$\Gamma(1-s+\nu)$  in the latter integral, we change the contour
$\sigma$ on the right-hand infinite loop $L_+$, encircling   the
right-hand simple poles of the obtained gamma functions in the
numerator.   Precisely, with a simple substitution we find
\begin{align}\label{eqn25}
&{1\over 4\pi i} \int_\sigma  { F(1-s) \  \Gamma(1-s+\nu)\over
\Gamma( (s-1/2)/2) \Gamma( (1/2-s)/2)}\xi^{-s} ds\nonumber \\&= -
{2^\nu \over 2\pi i \sqrt\pi } \int_{L_+} { F(1-2s) \
\Gamma((1+\nu)/2-s) \Gamma(1+\nu/2- s) \over \Gamma( s-1/4) \Gamma(
1/4-s)} \left(4 \xi^2\right)^{- s} ds.
\end{align}

Meanwhile, the integral
\begin{align}
-  {1 \over 2\pi i } \int_{L_+} { \Gamma((1+\nu)/2-s)
\Gamma(1+\nu/2- s) \over \Gamma( s-1/4) \Gamma( 1/4-s)} z^{- s} ds
\end{align}
can be calculated, using the Slater theorem \cite{prud}, Vol. 3, and
it can be expressed in terms of the hypergeometric functions
${}_1F_{2}$.  Consequently,  we obtain
\begin{align}\label{eq}
-&  {1 \over 2\pi i } \int_{L_+} {  \Gamma((1+\nu)/2-s)
\Gamma(1+\nu/2- s) \over \Gamma( s-1/4) \Gamma( 1/4-s)} z^{- s} ds =
{\sin(\frac{\pi}{2}(\nu+\frac32) ) \over \sqrt \pi z^{1+\nu/2}}\nonumber \\
&\times  \left({3\over 2}+\nu\right)   {}_1F_{2} \left({7\over
4}+{\nu\over 2}; {3\over 4} + {\nu\over 2}, {3\over 2}; -{1\over z}
\right)- {\sin(\frac{\pi}{2}(\nu+\frac12) ) \over 2\sqrt
\pi z^{(1+\nu)/2}}\nonumber\\
&\times\left({1\over 2}+\nu\right) { }_1F_{2} \left({5\over
4}+{\nu\over 2};  {1\over 4} + {\nu\over 2}, {1\over 2}; -{1\over z}
\right).
\end{align}

However,  the right-hand side of \eqnref{eq}  can be simplified and
written in terms of the elementary trigonometric functions,
appealing to relation (7.14.1.1)  in \cite{prud}, Vol. 3 and keeping
in mind particular cases of the modified Bessel function of the
first kind.  Thus we, finally,  derive
\begin{align}
&-  {1 \over 2\pi i } \int_{L_+} {  \Gamma((1+\nu)/2-s)
\Gamma(1+\nu/2- s) \over \Gamma( s-1/4) \Gamma( 1/4-s)} z^{- s} ds =
\nonumber \\& {z^{-(1+\nu)/2}\over 2\sqrt \pi}\left[
\cos\left({2\over \sqrt z}\right) \left[ \left({1\over 2}+\nu\right)
\sin\left({\pi\over 2} \left(\nu+{1\over 2}\right)\right)+ {\cos
(\pi\nu/2)\over \sqrt z} \right]\right. \nonumber \\  & \left.  -
\sin \left({2\over \sqrt z}\right) \left[ \left({3\over
2}+\nu\right) \cos\left({\pi\over 2} \left(\nu+{1\over
2}\right)\right) +  {1\over 2} \cos \left({\pi\nu\over 2}\right)
\right.\right.\nonumber \\  &\left.\left.\quad + {1\over \sqrt z}
\sin \left({\pi\over 2} \left (\nu+ {1\over 2}\right)\right)\
\right] \right] .
\end{align}

Hence, combining with \eqnref{eqn21}, \eqnref{eqn24},
\eqnref{eqn25}, after straightforward manipulations we arrive at the
final form \eqnref{eqn10} of solutions for the integral equation
\eqnref{eqn1}, completing the proof of Theorem 1.

\end{proof}

Some simplification may be made to our solution.

\begin{cor}
When  ${\rm Re}(\nu) > -1/2$,  $\nu\neq 2n+1/2, \ n \in
\mathbb{N}_0$, the solutions  \eqnref{eqn10} take the form
\begin{align}\label{eqn26}
&\varphi(\xi) = - {1\over 4\pi i} \int_\sigma {\rho(s)\over
\Gamma(1-s+\nu) }\ \xi^{-s} ds +
\frac{\Gamma(\frac12-\nu)}{\xi^{(1+\nu)}}
\sin\left(\frac{\pi}{2}(\nu+\frac12) \right)\times\nonumber \\&
{1\over \pi\sqrt \pi } \int_{-\infty}^\infty  e^{(\nu+\frac12)u}
\left[ \left({1\over \sqrt 2}+\sqrt 2\ \nu\right)
\sin\left({\pi\over 2} \left(\nu+{1\over 2} \right) -
\frac{e^{u}}{\xi}\right) +  \cos
\left(\frac{\pi\nu}{2}\right)\right.\nonumber\\
& \times \left[  {  e^{u} \over \sqrt 2 \xi}\ \left(
\cos\left(\frac{e^{u}}{\xi} \right)  -  {1\over \sqrt 2} \sin
\left(\frac{e^{u}}{\xi} \right) \right)  -  \sin
\left(\frac{e^{u}}{\xi} \right)  \left(1+ {1\over  \sqrt 2}\right)
\right]\nonumber \\ &\left.\quad\quad+  \sin \left({\pi\nu \over 2}
\right) \sin \left(e^{u} \xi^{-1} \right) \left( 1- {e^{u} \over  2
\xi}\right)\right]\times\nonumber \\ &   \int_0^\infty \left[
\left(\sinh u+\sinh x\right)^{\nu-1/2} + \left|\sinh u-\sinh x
\right |^{\nu-1/2}\right]\frac{ f(x)}{ (\cosh x)^{\nu-1}}  \ dx du,
\end{align}
for $\xi>0.$
\end{cor}
\begin{proof}
The proof will be completed if we compute the inner integral with
respect to $y$ in \eqnref{eqn10} after the corresponding interchange
of the order of integration. In fact, let $|{\rm Re}(\nu) | < 1/2$.
Then since $ e^{3x/2} \  f (x) \in L_1(\mathbb{R}_+)$, the integral
with respect to $x$ in \eqnref{eqn10} converges absolutely and
uniformly for $y \ge 0$. Hence the interchange is possible, where
the relatively convergent integral by $y$ can be calculated via
relation (2.5.3.10) in \cite{prud}, Vol. 1.  This leads to
\eqnref{eqn26}. Moreover, the obtained integral with respect to $x$
in the right-hand side of \eqnref{eqn26} converges absolutely and
uniformly for $\nu,\ {\rm Re}(\nu) \ge \nu_0
> - 1/2$ under the condition $ e^{3x/2}   f (x)  \in
L_1(\mathbb{R}_+)$, representing an  analytic function of $\nu$ in
the domain $D= \{ \nu \in \mathbb{C}:  {\rm Re}( \nu) > - 1/2,\  \nu
\neq 2n+ 1/2,\ n \in \mathbb{N}_0\}$.   Hence,   by analytic
continuation   equality \eqnref{eqn26} holds for all $\nu \in D$.
\end{proof}

\subsection{Conditions for Uniqueness}
We now turn to the question of uniqueness of solutions for the
integral equation. \eqnref{eqn1}. One way to obtain this is to
assume the the Mellin transform $\varphi^*(s)$ satisfies a
particular condition on $\sigma.$

\begin{cor}\label{unique}

Suppose that $f$ and $\varphi$ satisfies the conditions of Theorem
\ref{Main} and the Mellin transform of $\varphi$ has the additional
property that
\begin{align}\label{product}
\varphi^*(s) \Gamma (1-s+\nu)=  \varphi^*(1-s) \Gamma (s+\nu),\ s
\in \sigma
\end{align}
i.e. this product is even on $\sigma$. Then the integral equation
\eqnref{eqn1} has a unique solution
\begin{align}
\varphi(\xi)=\frac{1}{2\pi i}\int_\sigma
\frac{F(s)\Gamma(s+\nu)\Gamma(1-s+\nu)}{2\Gamma\left(\frac{1/2-s}{2}\right)\Gamma\left(\frac{s-1/2}{2}\right)}\xi^{-s}ds,
\end{align}
where $F$ is given by \eqnref{eqn21}.
\end{cor}
\begin{rem}
Equality \eqnref{product} can be rewritten, for instance,  in terms
of the modified Laplace transforms  of the function $\varphi$.
Indeed, recalling the Mellin-Parseval identity \eqnref{eqn7},
equality \eqnref{product} can be written in the  following
equivalent form
$$ x^{-1-\nu} \int_0^\infty \exp( - t/x) \varphi(t)  t^\nu  dt=  x^\nu \int_0^\infty \exp( - xt) \varphi(t)  t^\nu dt,\ x >0.$$
\end{rem}

\begin{proof}
Using the condition on the Mellin transform $\varphi^*$, we find
from equation \eqnref{eqn532} that
\begin{align}
\varphi^*(s)&=\frac{[F(1-s)+F(s)]\Gamma(s+\nu)\Gamma(1-s+\nu)}{4\Gamma\left(\frac{1/2-s}{2}\right)\Gamma\left(\frac{s-1/2}{2}\right)}\nonumber
\\&=\frac{F(s)\Gamma(s+\nu)\Gamma(1-s+\nu)}{2\Gamma\left(\frac{1/2-s}{2}\right)\Gamma\left(\frac{s-1/2}{2}\right)},
\end{align}
since $F(s)=F(1-s)$ on $\sigma$. Inverting the Mellin transform we
end up with the unique solution in the form,
\begin{align}
\varphi(\xi)=\frac{1}{2\pi i}\int_\sigma
\frac{F(s)\Gamma(s+\nu)\Gamma(1-s+\nu)}{2\Gamma\left(\frac{1/2-s}{2}\right)\Gamma\left(\frac{s-1/2}{2}\right)}\xi^{-s}ds.
\end{align}

\end{proof}

A natural question is when $\psi^*(s)=\varphi^*(s) \Gamma
(1-s+\nu)=\psi^*(1-s)$ on $\sigma?$

\begin{lem} Let $\psi$ be the
inverse Mellin transform of $$\psi^*(s)=\varphi^*(s)\Gamma
(1-s+\nu)=\psi^*(1-s).$$ If $g(y)=\psi(e^y)e^{1/2y}$ is even, then
$\psi^*(s)=\psi^*(1-s)$ for $s\in \sigma.$
\end{lem}
\begin{proof}
If $s\in \sigma$, then $s=\frac12+i\tau.$ Clearly we require
\begin{align}
\int_0^\infty x^{-1/2+i\tau}\psi(x)dx=\int_0^\infty
x^{-1/2-i\tau}\psi(x)dx,
\end{align}
for every $\tau\in \Real.$ We convert this to a Fourier transform by
setting $x=e^y$. This gives
\begin{align}
\int_{-\infty}^\infty e^{i\tau
y}\psi(e^y)e^{1/2y}dy=\int_{-\infty}^\infty e^{-i\tau
y}\psi(e^y)e^{1/2y}dy.
\end{align}
So the Fourier transform of $g(y)=\psi(e^y)e^{1/2y}$ is even. This
happens precisely when $g$ is even.
\end{proof}

We are then able to write down a solution of our Cauchy problem.
\begin{thm}
Let ${\rm Re}(\nu) > -1/2, \varphi(\xi) \in  S (\mathbb{R}_+),\
e^{3x/2} \  f (x)  \in L_1(\mathbb{R}_+), \ $ and $y^{-\nu-1/2} (F_c
g)( y) \in \mathcal{S} (\mathbb{R}_+)$,\ where $ g(x)=
(x^2+1)^{-\nu/2} \ f \left(\sinh^{-1}  x\right)$. Let also the
Mellin transform $\varphi^*$ of $\varphi$ satisfy the conditions

$$\varphi^*(s)  |s|^{-\nu}  e^{(\pi/2 -\delta) |s|} \in L_1(\sigma),\  \delta \in [0, \pi/2) ,$$
and
\begin{align}
\varphi^*(s) \Gamma (1-s+\nu)=  \varphi^*(1-s) \Gamma (s+\nu),\ s
\in \sigma
\end{align}
where  $\sigma = \left\{ s\in \mathbb{C}, \  s = {1\over 2}+i\tau,\
\tau \in \mathbb{R}\right\}$ and
\begin{align}\label{eqn9}
h_{c}(\tau)&=   \varphi^*(1/2 +i\tau)  \Gamma(1/2- i\tau+\nu)  +
\varphi^*(1/2 -i\tau) \Gamma(1/2+ i\tau+\nu) \\&=
O(\tau^2),\nonumber
\end{align}
$\tau \to 0.$ Define
\begin{align}\label{eqn21}
F(s) &=     \int_1^\infty  \left[ \zeta^{s-\frac32} +
\frac{1}{\zeta^{s+\frac12}} \right] \int_0^\infty \frac{  \cos
(\frac12\eta(\zeta-\frac{1}{\zeta}) )}{ \eta^{1/2+\nu}} \nonumber
\\&\times \int_0^\infty \frac{\cos( \eta\sinh \xi)
}{(\cosh \xi)^{ \nu-1/2}}f(\xi) d\xi d\eta d\zeta.
\end{align}
Then the problem
\begin{align*}
u_{tt}&=u_{xx}+(\mathrm{sech}^2x)u,\\
u(x,0)&=f(x),\ u_t(x,0)=0,
\end{align*}
has a solution
\begin{align}\label{finintegral}
u(x,t)&=\frac{1}{ \pi i}\int_0^\infty  \int_\sigma
K(t,x,\xi)\frac{F(s)\Gamma(s+\nu)\Gamma(1-s+\nu)}{\Gamma\left(\frac{1/2-s}{2}\right)\Gamma\left(\frac{s-1/2}{2}\right)}\xi^{-s}ds
d\xi,
\end{align}
where
\begin{align}
K(t,x,\xi)&= \frac{ \cos(e^{t} \xi \sinh
    x)e^{ \frac{t}{2}} J_{\nu +\frac{1}{2}}\left(e^{t} \xi
   \cosh  x \right)}{\sqrt{\xi  \mathrm{sech} x }}\nonumber  \\& \quad +\frac{ \cos(e^{-t} \xi  \sinh
    x)e^{ -\frac{t}{2}} J_{\nu +\frac{1}{2}}\left(e^{-t} \xi
   \cosh  x \right)}{\sqrt{\xi  \mathrm{sech} x }}.
\end{align}
\end{thm}
\begin{proof}
We only need to establish convergence of the final integral and this
follows from the conditions of the theorem. Specifically, since the
Fourier cosine transform maps Schwartz functions to Schwartz
functions, and the Mellin transform in $s$ of a Schwartz function,
with $\mathrm{Im}(s)=t$, is Schwartz in $t$, the function $F$ has
rapid decay on the line $\sigma$ and hence its inverse Mellin
transform has rapid decay. This guarantees the convergence of the
final integral.
\end{proof}

Simplification of this result may be possible, but we will not
discuss this question here.

\subsection{The Second Integral Equation}

We finish by turning now to the solution of the integral equation
\eqnref{eqn2}. This can be established in the same manner as was
used for the proof of Theorem \ref{Main} with the use of the
integral \eqnref{eqn4} and is given by our next result.

\begin{thm}\label{Main2}
Let ${\rm Re}(\nu) > -\frac12, \varphi(\xi) \in
\mathcal{S}(\mathbb{R}_+),\ e^{3x/2} \  f (x)  \in
L_1(\mathbb{R}_+), \  $ and  $y^{-\nu-1/2} (F_s  g)( y) \in
\mathcal{S} (\mathbb{R}_+)$,\ where $ g(x)= (x^2+1)^{-\nu/2} \ f
\left(\sinh^{-1}(\  x)\right)$. Let also the Mellin transform
$\varphi^*$ of $\varphi$ satisfy the condition
\begin{align}
\varphi^*(s)  |s|^{-\nu}  e^{(\pi/2 -\delta) |s|} \in L_1(\sigma),\
\delta \in [0, \pi/2) ,
\end{align}
  where  $\sigma = \left\{ s\in \mathbb{C}, \  s = {1\over 2}+i\tau,\
\tau \in \mathbb{R}\right\}$ and

\begin{align}\label{eqn27}
h_{s}(\tau)&=   \varphi^*(1/2 +i\tau)  \Gamma(3/2- i\tau+\nu)  +
\varphi^*(1/2 -i\tau) \Gamma(3/2+ i\tau+\nu)\nonumber \\
& = O(\tau^2),\   \tau \to 0.
\end{align}

Then the integral equation \eqnref{eqn2} has the  following
solutions
\begin{align}\label{eqn28}
\varphi(\xi) & = - {1\over 4\pi i} \int_\sigma {\rho(s)\over
\Gamma(2-s+\nu) }\ \xi^{-s} ds
 + {\sqrt 2 \  \xi^{-(2+\nu)} \over \pi\sqrt \pi }\nonumber \\ &\times \int_{-\infty}^\infty  e^{(\nu+3/2)u}  \sin\left({\pi\over 2} \left(\nu+{1\over 2} \right) - \frac{e^{u}}{  \xi} \right) \cosh u
\nonumber \\
&\times  \int_0^\infty \   y^{1/2-\nu} \cos (y\sinh u) \int_0^\infty
\sin ( y\sinh x)  (\cosh x)^{1-\nu} f(x) \ dx dy du,
\end{align}
for $\xi>0$,  depending on an arbitrary function  $\rho(s)$,  which
satisfies relation \eqnref{eqn11} on $\sigma$.
\end{thm}

\begin{proof}  In fact, an analog of equation \eqnref{eqn13} will be the equality
\begin{align}\label{eqn29}
&\sqrt{\pi\over 2}  {\sinh x \ \cosh^{\nu} x \over 2\pi i}
\int_{\sigma} \varphi^*(s)   \frac{ 2 \Gamma\left(\frac{ 2-s+\nu
}{2}\right)}{2^s\Gamma\left(\frac{ 2+\nu+s }{2}\right)}\times
\nonumber\\&{}_2F_1 \left( {2-s+\nu\over 2}, {1+\nu + s\over 2};
{3\over 2};  - \sinh^2 x \right)\    ds = f(x).
\end{align}
Hence,  by a similar calculation to that leading  to \eqnref{eqn14},
we get
\begin{align}\label{eqn30}
 \int_{0}^\infty  h_{s}(\tau)\frac{    {}_2F_1 \left( {3/2+\nu+i\tau \over 2},
{3/2+\nu -i\tau \over 2}; {3\over 2};  - \sinh^2 x \right)}{\left[
\Gamma\left({5/2+ \nu+i\tau\over 2}\right) \Gamma\left({5/2+ \nu-
i\tau\over 2}\right)\right]}\  d\tau &= 2^{3/2+\nu}  {\cosh^{-\nu}
x\over \sinh x}  \ f(x)
\end{align}
and since
\begin{align}
&\sinh x \ {}_2F_1 \left( {3/2+\nu+i\tau \over 2}, {3/2+\nu -i\tau
\over 2}; {3\over 2};  - \sinh^2 x \right) \nonumber =\\&
 \frac{2^{1/2-\nu}  }{
\Gamma((3/2+ \nu+i\tau)/2) \Gamma((3/2+ \nu- i\tau)/2)}\int_0^\infty
y^{\nu-1/2} \sin (y\sinh x) K_{i\tau} (y) dy,
\end{align}
${\rm Re}( \nu) > -1/2$, we derive an analog of the equation
\eqnref{eqn16}
\begin{align}\label{eqn31}
\int_0^\infty \int_{0}^\infty    {y^{\nu -1/2}h_{s}(\tau) \sin(u y)
\ K_{i\tau} (y) \over \Gamma(\frac32+ \nu+i\tau) \Gamma(\frac32+
\nu- i\tau)} d\tau dy = \pi (u^2+1)^{-\nu/2}  \  f \left(\sinh^{-1}
( u)\right).
\end{align}

Further, according to conditions of the theorem we may take the
Fourier sine transform and then the  Kontorovich-Lebedev transform,
to the function $h_s(\tau)$ in the form
\begin{align}\label{eqn32}
 h_{s}(\tau)& =  {4\over \pi^2} \tau\sinh(\pi\tau) \left|\Gamma\left(\frac32+ \nu+i\tau\right) \right|^2  \int_0^\infty  y^{-1/2-\nu} K_{i\tau}
 (y)\nonumber
\\
& \times   \int_0^\infty  \sin ( y\sinh x)  (\cosh x)^{1-\nu} f(x) \
dx dy.
\end{align}
Moreover, after the use of the  Parseval equality for the Fourier
sine transform \cite{tit},  relation (2.16.14.1) in \cite{prud},
Vol. 2, integration by parts and differentiation under the integral
sign in the integral with respect to $u$,  which is permitted under
conditions of the theorem, we obtain
\begin{align}\label{eqn33}
h_{s}(\tau)& =  {4\sqrt 2 \over \pi\sqrt\pi}  \cosh
\left({\pi\tau\over 2}\right) \left|\Gamma\left(\frac32+
\nu+i\tau\right)\right|^2 \int_0^\infty\int_0^\infty   \int_0^\infty
\cos (\tau u) \cosh u \nonumber \\& \times    \ y^{1/2-\nu} \cos
(y\sinh u)\sin ( y\sinh x)  (\cosh x)^{1-\nu} f(x) \ dx dy du.
\end{align}

Then similarly to equations \eqnref{eqn20}, \eqnref{eqn21} and
\eqnref{eqn22}  we derive from \eqnref{eqn33}

\begin{align}\label{eqn34}
 \varphi^*(s) {\Gamma(2-s+\nu)\over \Gamma(1+ s+\nu)}  +
\varphi^*(1-s) = \frac{ G(s)\  \Gamma( 2-s+\nu)}{ \Gamma( (s+1/2)/2)
\Gamma( (3/2-s)/2)}, \  s \in \sigma,
\end{align}
where
\begin{align}\label{eqn35}
G(s)=  {2\sqrt 2 \over \sqrt\pi}   \int_{-\infty}
^\infty\int_0^\infty\int_0^\infty & e^{u(1/2-s)}   \cosh u  \
y^{1/2-\nu} \cos (y\sinh u)\nonumber \\& \times      \sin( y\sinh x)
(\cosh x)^{1-\nu} f(x) \ dx dy du,
\end{align}
Hence as in the proof of Theorem \ref{Main},    the final solutions
of the equation \eqnref{eqn2} can be expressed in terms of the
inverse Mellin transform  \eqnref{eqn6}
\begin{align}\label{eqn36}
\varphi(\xi) = - {1\over 4\pi i} \int_\sigma {\rho(s)\xi^{-s}\over
\Gamma(2-s+\nu) }\  ds + {1\over 4\pi i} \int_\sigma { G(1-s) \
\Gamma(2-s+\nu)\xi^{-s}\over \Gamma\left(\frac{ s+\frac12}{2}\right)
\Gamma\left(\frac{  \frac32-s }{2}\right)} ds.
\end{align}
Then as above,
\begin{align}\label{eqn37}
{1\over 4\pi i}& \int_\sigma  { G(1-s) \ \Gamma(2-s+\nu)\over
\Gamma( (s+1/2)/2) \Gamma( (3/2-s)/2)}\xi^{-s} ds =\nonumber \\&  -
{2^{\nu+1} \over 2\pi i \sqrt\pi } \int_{L_+} { G(1-2s) \  \Gamma(1
+\nu/2-s) \Gamma(3/2+\nu/2- s) \over \Gamma( s+1/4) \Gamma( 3/4-s)}
\left(4 \xi^2\right)^{- s} ds.
\end{align}

Meanwhile, the integral
$$-  {1 \over 2\pi i } \int_{L_+} {  \Gamma(1+\nu/2-s) \Gamma(3/2+\nu/2- s) \over \Gamma( s+1/4) \Gamma( 3/4-s)} z^{- s} ds$$
can be calculated in terms of the elementary functions.  Precisely,
we derive
\begin{align}
-  {1 \over 2\pi i } \int_{L_+} {  \Gamma(1+\frac{\nu}{2}-s)
\Gamma(\frac32+\frac{\nu}{2}- s) \over \Gamma( s+\frac14)
\Gamma(\frac34-s)} z^{- s} ds =  \frac{\sin\left({\pi\over
2}\left(\nu + {\frac12}\right) - {2\over \sqrt z}
\right)}{\sqrt{\pi} z^{1+\frac{\nu}{2}}}.
\end{align}
Hence, combining with \eqnref{eqn35}, \eqnref{eqn36} and
\eqnref{eqn37}, we arrive at the final form \eqnref{eqn28} of
solutions for the integral equation \eqnref{eqn2}, completing the
proof of Theorem \ref{Main2}.

\end{proof}

The corresponding corollary for the values $\nu \in \mathbb{C} : \
{\rm Re}(\nu) > -1/2,  \  \nu\neq  2n+7/2, \ n \in \mathbb{N}_0$
can be formulated as follows

\begin{cor} When  ${\rm Re}(\nu) > -\frac12,\nu\neq 2n+7/2,
\ n \in \mathbb{N}_0$, the solutions \eqnref{eqn28} take the form
$$\varphi(\xi) = - {1\over 4\pi i} \int_\sigma {\rho(s)\over   \Gamma(2-s+\nu) }\ \xi^{-s} ds +   \xi^{-(2+\nu)} \Gamma(3/2-\nu) \sin(\pi(3/2-\nu)/2)  $$
$$ \times  { 1 \over \pi\sqrt {2\pi} } \int_{-\infty}^\infty  e^{(\nu+3/2)u}  \sin\left({\pi\over 2} \left(\nu+{1\over 2} \right) - e^{u}  \xi^{-1}\right) \cosh u $$
$$\times \left[ \int_0^\infty  \left( \sinh u + \sinh x\right)^{\nu-3/2}  (\cosh x)^{1-\nu} f(x) \ dx  du\right.$$
$$\left. -   \int_0^u  \left( \sinh u - \sinh x\right)^{\nu-3/2}  (\cosh x)^{1-\nu} f(x) \ dx  du \right.$$
$$\left. +  \int_u^\infty  \left( \sinh x - \sinh u \right)^{\nu-3/2}  (\cosh x)^{1-\nu} f(x) \ dx  du\right], \  \xi >0.$$

\end{cor}

We conclude with a condition that will guarantee unique solutions
for the second integral equation \eqnref{eqn2}. The proof is similar
to that for Corollary \ref{unique} and we omit it.

\begin{cor}
Let $f$ and $\varphi$ satisfy the conditions of Theorem \ref{Main2}
and suppose further that

$$\varphi^*(s) \Gamma (2-s+\nu)=  \varphi^*(1-s) \Gamma (1+ s+\nu),\
s \in \sigma.$$ Then equation \eqnref{eqn2} has a unique solution
given by
\begin{align}
\varphi(\xi)=\frac{1}{4\pi i}\int_\sigma
\frac{G(1-s)\Gamma(1+s+\nu)\Gamma
(2-s+\nu)}{\Gamma\left(\frac{3/2-s}{2}\right)\Gamma\left(\frac{s+1/2}{2}\right)}\xi^{-s}ds,
\end{align}
where $G$ is given by \eqnref{eqn35}.
\end{cor}

\section{Conclusion}
The methods introduced in \cite{Cra2014} and extended in the current
work lead to many interesting problems involving integral transforms
and integral equations. There are still very many open questions,
even in terms of the focus of this paper. For example, one would
like to be able to solve the general integral equation
\eqnref{generalinteqn}. When applied to other PDEs, the method also
generates integral equations that do not seem to have been studied
in the literature and many interesting questions arise from the
investigation of these problems. We hope that this work stimulates
further research in this area.

\section*{Acknowledgment}
The second  author was partially supported by CMUP (UID/MAT/
00144/2013), which is funded by FCT (Portugal) with national (MEC)
and European structural funds through the programs FEDER, under the
partnership agreement PT2020. The first author would like to thank
Professor Peter Olver for his hospitality at the University of
Minnesota in September of 2013, and for reminding the author that
one can subtract as well as add.

\bibliographystyle{plain}
\bibliography{C:/Mark/Transfer/TexFiles/include/bibliography}
\end{document}